\documentclass[11pt,draftcls,onecolumn]{IEEEtran}

\pdfoutput=1  

\IEEEoverridecommandlockouts                              
\usepackage{cite,url}

\usepackage[usenames,dvipsnames]{color}
\usepackage{graphicx,wrapfig,subfigure}

\usepackage[cmex10]{amsmath}
\usepackage{amsfonts,amssymb,mathrsfs}
\usepackage[colorlinks,bookmarksopen,bookmarksnumbered,citecolor=red,urlcolor=red]{hyperref}

\newlength{\noteWidth}
\long\def\notes#1{\ifinner
           {\footnotesize #1}
           \else
           \marginpar{\parbox[t]{\noteWidth}{\raggedright\footnotesize #1}}
       \fi\typeout{#1}}

\def\notes#1{\typeout{read notes: #1}}  

\def\spm#1{\notes{SPM:  #1}}

\def\sq{\hbox{\rlap{$\sqcap$}$\sqcup$}}
\def\qed{\ifmmode\sq\else{\unskip\nobreak\hfil
\penalty50\hskip1em\null\nobreak\hfil\sq
\parfillskip=0pt\finalhyphendemerits=0\endgraf}\fi\medskip}

\long\def\defbox#1{\framebox[.9\hsize][c]{\parbox{.85\hsize}{%
\parindent=0pt
\baselineskip=12pt plus .1pt      
\parskip=6pt plus 1.5pt minus 1pt 
 #1}}}

\long\def\beginbox#1\endbox{\subsection*{}%
\hbox{\hspace{.05\hsize}\defbox{\medskip#1\bigskip}}%
\subsection*{}}

\def\endbox{}

\def\transpose{{\hbox{\it\tiny T}}}

\newsavebox{\junk}
\savebox{\junk}[1.6mm]{\hbox{$|\!|\!|$}}

\def\V{{\sf V}}

\def\state{{\sf X}}

\newcommand{\field}[1]{\mathbb{#1}}

\def\Re{\field{R}}

\def\ind{\field{I}}

\def\cpi{\check{\pi}}

\def\cP{{\check{P}}}

\def\bfmath#1{{\mathchoice{\mbox{\boldmath$#1$}}%
{\mbox{\boldmath$#1$}}%
{\mbox{\boldmath$\scriptstyle#1$}}%
{\mbox{\boldmath$\scriptscriptstyle#1$}}}}

\def\bfmr{\bfmath{r}}

\def\bfmy{\bfmath{y}}

\def\bfmX{\bfmath{X}}

\def\bfmtily{\widetilde{\bfmath{y}}}

\def\tily{\tilde{y}}

\def\bfmY{\bfmath{Y}}

\def\bfmhhaY{\bfmath{\hhaY}} 
\def\bfmhhaY{\hbox to 0pt{$\widehat{\bfmY}$\hss}\widehat{\phantom{\raise 1.25pt\hbox{$\bfmY$}}}}

\def\bfatom{\bfmath{\theta}}
\def\bfPhi{\bfmath{\Phi}}

\def\bfzeta{\bfmath{\zeta}}

\def\haP{{\widehat P}}

\def\til={{\widetilde =}}

\def\clE{{\cal E}}
\def\clF{{\cal F}}
\def\clG{{\cal G}}
\def\clH{{\cal H}}

\def\clP{{\cal P}}

\def\bfatom{{\mbox{\boldmath$\alpha$}}}
\def\atom{{\mathchoice{\bfatom}{\bfatom}{\alpha}{\alpha}}}

\def\head#1{\subsubsection*{#1}}

 \def\FRAC#1#2#3{\genfrac{}{}{}{#1}{#2}{#3}}

\def\ddtp{{\mathchoice{\FRAC{1}{d^{\hbox to 2pt{\rm\tiny +\hss}}}{dt}}%
{\FRAC{1}{d^{\hbox to 2pt{\rm\tiny +\hss}}}{dt}}%
{\FRAC{3}{d^{\hbox to 2pt{\rm\tiny +\hss}}}{dt}}%
{\FRAC{3}{d^{\hbox to 2pt{\rm\tiny +\hss}}}{dt}}}}

\def\half{{\mathchoice{\FRAC{1}{1}{2}}%
{\FRAC{1}{1}{2}}%
{\FRAC{3}{1}{2}}%
{\FRAC{3}{1}{2}}}}

\def\eqdef{\mathbin{:=}}

\def\Prob{{\sf P}}

\def\Expect{{\sf E}}

\def\average#1,#2,{{1\over #2} \sum_{#1}^{#2}}

\def\eye(#1){{\bf(#1)}\quad}

\newtheorem{theorem}{Theorem}[section]

\newtheorem{proposition}[theorem]{Proposition}
\newtheorem{lemma}[theorem]{Lemma}

\def\Lemma#1{Lemma~\ref{#1}}
\def\Proposition#1{Proposition~\ref{#1}}
\def\Theorem#1{Theorem~\ref{#1}}

\def\Section#1{Section~\ref{#1}}

\def\eq#1/{(\ref{e:#1})}

\newcommand{\beqn}[1]{\notes{#1}%
\begin{eqnarray} \elabel{#1}}

\newcommand{\eeqn}{\end{eqnarray} }

\newcommand{\beq}[1]{\notes{#1}%
\begin{equation}\elabel{#1}}

\newcommand{\eeq}{\end{equation}}

\def\bdes{\begin{description}}
\def\edes{\end{description}}

\def\barP{{\overline {P}}}

\newcounter{rmnum}
\newenvironment{romannum}{\begin{list}{{\upshape (\roman{rmnum})}}{\usecounter{rmnum}
\setlength{\leftmargin}{14pt}
\setlength{\rightmargin}{8pt}
\setlength{\itemsep}{2pt}
\setlength{\itemindent}{-1pt}
}}{\end{list}}

\newcounter{anum}

\def\Spx{\textsf{S}}

\def\bfgamma{\bfmath{\gamma}}

\def\bfgamma{\bfmath{\gamma}}

\def\welf{\mathchoice{\mbox{\small$\cal W$}}%
{\mbox{\small$\cal W$}}%
{\mbox{$\scriptstyle\cal W$}}%
{\mbox{$\scriptscriptstyle\cal W$}}}

\def\util{\mathchoice{\mbox{\small$\cal U$}}%
{\mbox{\small$\cal U$}}%
{\mbox{$\scriptstyle\cal U$}}%
{\mbox{$\scriptscriptstyle\cal U$}}}

\def\tilutil{\mathchoice{\mbox{\small$\cal \widetilde U$}}%
{\mbox{\small$\cal\widetilde U$}}%
{\mbox{$\scriptstyle\cal \widetilde U$}}%
{\mbox{$\scriptscriptstyle\cal \tilde U$}}}

\def\Ebox#1#2{%
\begin{center}
\includegraphics[width= #1\hsize]{#2} \end{center}}

\def\varH{{\cal S}}
\def\derH{H}

\def\avar{\kappa^2}

\def\cp{\check p}

\def\DV{{\cal H}}
\def\V{{\cal V}}
\def\GNoplus{g_{\text{\scriptsize p}}} 

\def\Fig#1{Fig.~\ref{#1}}

\def\ind{\field{I}}

\def\Re{\field{R}}

\title{Ancillary Service to the Grid\\
Using Intelligent Deferrable Loads
}

\author{Sean Meyn, Prabir Barooah,  
Ana Bu\v{s}i\'{c},
Yue Chen,
and Jordan Ehren
\thanks{This research is supported by the NSF grant CPS-0931416, the Department of Energy Awards DE-OE0000097 \&\ DE-SC0003879,
and the French National Research Agency grant ANR-12-MONU-0019.  We acknowledge the help of Mark Rosenberg
who offered many suggestions to improve the manuscript, and caught several typos in earlier drafts.}
\thanks{S.M.,  Y.C.\  and J.E.\ are with the Dept.\ ECE and P.B. is with the Dept.\ MAE at the University of Florida, Gainesville.   
A.B.\ is with  INRIA and the Computer Science Dept. of \'Ecole Normale Sup\'erieure, Paris, France. 
}
}

\begin{document}

\maketitle
\thispagestyle{empty}

\begin{abstract} 
Renewable energy sources such as wind and solar power have a high
degree of unpredictability and time-variation, which makes balancing
demand and supply challenging. 
One possible way to address this challenge is to
harness the inherent flexibility in demand of many types of loads.  
Introduced in this paper is a technique for decentralized control
for automated demand response that can be used by grid operators
as ancillary service for maintaining demand-supply balance.

A  Markovian Decision Process (MDP) model is introduced for  an
individual load.  A randomized control architecture is proposed,
motivated by the need for decentralized decision making, and the need
to avoid synchronization that can lead to large and detrimental spikes
in demand.  An aggregate model for a large number of loads is then
developed by examining the mean field limit.  A key innovation is an
LTI-system approximation of the aggregate nonlinear model, with a
scalar signal as the input and a measure of the aggregate demand as
the output. This makes the approximation particularly convenient for control design at the grid level.
 
The second half of the paper contains a detailed    application of these results to a network of residential pools. 
Simulations are provided to illustrate the accuracy of the approximations and effectiveness of the proposed control approach.


\end{abstract}

%

\clearpage

\section{Introduction} 
\label{s:intro}

Renewable energy penetration is rising rapidly throughout the world,
and bringing with it high volatility in energy supply.   Resources
are needed to compensate for these large fluctuations in power.  
The federal energy regulatory commission (FERC)
in conjunction with generation and utility companies are struggling to
find resources, and finding ways to properly compensate for ancillary
services that are badly needed by each \textit{balancing authority} (BA) in the U.S..  
FERC orders 755 and 745 are examples of their attempts to
provide incentives. 
\notes{the utilities are often victims in all of this -- I hope BA is o.k.  -spm\\
Expand on FERC.   Also see refs from Toulouse (in commented text)
} 

%

This paper concerns decentralized control of a large number of electric loads in a power grid.
A particular load has a service it is intended to provide -- clean dishes, hot water, or a clean pool.
It is assumed that each load has some flexibility in energy consumption.  This flexibility is harnessed to
 provide ancillary services to the power grid to help maintain stability, and to help offset any volatility in the grid because of line or generation outage, or because of the volatile nature of renewable energy.  This is commonly called ``demand response'', but the meaning is slightly different here:  The tuning of energy consumption is automated, and we assume that the consumers do not suffer any degradation in the service offered by the loads.

We argue that most of the load
in the U.S. is highly flexible, and this flexibility can be harnessed
to provide ancillary service without central control, and without
significant impact on the needs of consumers or industry. 
A defining characteristic of ancillary service is that on average it is a
\emph{zero-energy} service, so that the desired power
consumption level to be tracked is zero on average.
This makes use of  deferrable loads   particularly attractive as sources of ancillary service.

Many utilities already employ demand response programs that use
deferrable loads to reduce peak demand and manage emergency situations. Florida Power and Light  (FPL), for example, has 780,000 customers enrolled in their \textit{OnCall Savings Program} in which residential air conditioners, water heaters, and pool pumps systems are automatically controlled when needed \cite{FPLsaving}. 
Today,  FPL  uses this service only 3--4 times per year \cite{FPLsaving}. While a valuable service to the grid, there is tremendous additional potential from these sources that today is virtually untapped.

Nearly all of America's ISOs/RTOs also
allow for demand side resources to participate in their regulation and spinning reserve markets, but as of the summer of 2013, 
only PJM    allows aggregation (with approval) \cite{maccapcalkil12}. Growth of these resources in these wholesale markets has helped lower costs per megawatt-hour from 2009 to 2011 \cite{maccapcalkil12}.  Still, markets for regulation and spinning reserves from traditional generation sources
continue to grow because of increasing dependency on renewable generation.

\begin{figure}[h]
\vspace{-.15cm}
\Ebox{.75}{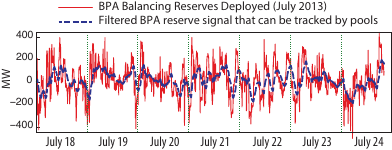} 
\vspace{-.25cm}
\caption{\textit{BPA Balancing Reserves Deployed}
---
Ancillary service needs at the BPA during one week in 2013.   The maximum is approximately one-tenth of maximum load in this region.  
}
\label{fig:BPA}
\vspace{-08pt}
\end{figure}

\Fig{fig:BPA} shows the regulation signal for a typical week within the Bonneville Power Authority (BPA)~\cite{BPA}.  
Its role 
is
analogous to the control signal in the feedback loop in a flight control system.   
Just like in an aviation control  system,  the variability seen in this figure is in part a consequence of  variability of wind generation in this region.  

We propose to break up a regulation signal into frequency bands for the purposes of ancillary services provisioning by various resources.   In prior work   it is shown how heating and ventilation systems in commercial buildings can provide service in the high frequency band, corresponding to periods ranging from 3 minutes to one hour \cite{haokowbarmey13,haobarmidmey12,linbarmey13}.   At the lowest frequencies,  an important resource will be flexible manufacturing.  An example is Alcoa, 
that today provides 70MW of service to MISO by providing control over their aluminum smelting operation in Indiana.
Alcoa's service is provided continuously, and provides significant revenue to Alcoa and even greater benefits to the region managed by MISO.

The technical content of the paper starts with a control architecture designed to address privacy concerns and communication constraints.   
It is assumed that an individual load can view a regulation signal, much as we can view BPA's regulation signal online today.

To provide ancillary service in a specified frequency band, we argue that it is essential to introduce randomization at each load. Among many benefits, randomization avoids synchronization, much like randomized congestion avoidance protocols in communication networks.   First deployed nearly fifty years ago, ALOHA may be the first distributed communication protocol based on randomization. 
\textit{Random Early Detection} for congestion control was introduced in the highly influential paper \cite{flojac93}. The historical discussion in this paper points to significant research on randomized algorithms beginning in the early 1970s, following wide deployment of ALOHA. Randomized protocols are now standard practice in communication networks \cite{sri04a}. It is likely that randomized algorithms will become a standard feature of the power grid of the future.

%
%
%
%


To formulate a randomized control strategy, a Markovian Decision Process (MDP) model is proposed for an individual load. 
An aggregate model for a large number of loads is then obtained as a mean field limit.

A particular formulation of  Todorov~\cite{tod07} is adopted because we can obtain an explicit solution,  and because of available tools 
for analysis borrowed from the theory of large deviations.   In particular, a key innovation in the present paper is an LTI--system approximation of the aggregate nonlinear model, which is possible through application of results from \cite{konmey05a}.    The scalar input in this linear model is a parameter that appears in the MDP cost function.

The LTI approximation is convenient for control design at the grid level:   the  input becomes the control signal that the BA will broadcast to all the loads, which adjusts a parameter in the randomized policy for the optimal MDP solution at each load.

In the second half of this paper we apply the general results of this paper to show how pool pumps can be harnessed to obtain  ancillary service in a medium frequency band, corresponding to the dashed line in \Fig{fig:BPA}. This is the same BPA regulation signal, passed through a low pass filter
 
 A pool pump is the
 heart of a pool's filtration system:  It runs each day for a period of time range from 4 to 24 hours, and consumes over  1~KW of power when in operation  
 \cite{PPDR08}.
The ability to control just half of
 Florida's pool pumps amounts to over 500~MW of power! Much of the control
 infrastructure is already in place~\cite{hallo06}.  Still, constraints
 and costs must be satisfied.  These include run-times per day and per
 week, the cost of startup and shut down, as well as the total energy
 consumption. Moreover, there are privacy concerns and
 related communication constraints.  Consequently, control algorithms
 must be distributed so that most of the required intelligence resides
 at individual pool pumps.  In this paper we focus on constraints
 related to run-times per day, which is critical for keeping the water
 in the pool clean.  Privacy and communication constraints will be addressed through
 the distributed control architecture. 

A number of recent works have explored the potential for flexible
loads for providing ancillary service. These include commercial
building  thermostatic loads to provide ancillary service in
the time-scale of a few minutes  (see \cite{matkoccal13} and refs.\ therein),
 electric
vehicle charging~\cite{macalhis10,tombou10,coupertemdeb12,matkoccal13}
that can provide ancillary service in the time scale of a few hours,
and our own recent work on harnessing ancillary
service from commercial building HVAC~\cite{haokowbarmey13,haobarmidmey12,linbarmey13}.

Mean-field games have been employed for analysis of aggregate loads in several recent papers 
\cite{macalhis10,coupertemdeb12}. 
See \cite{huacaimal07,borsun12,gasgauleb12} for more on general theory of mean-field technques.

The work of~\cite{matkoccal13} is most closely related to the 
present  paper, in that the authors also consider an aggregate model for a large collection of loads.  The natural state space model is bilinear, and converted to a linear model through division of the state.  The control architecture consists of a
centralized control signal computation based on state feedback,
and the resulting input is broadcast to the devices.   

In this paper, intelligence is concentrated at the individual load:  An
MDP control solution is obtained at each load, but the aggregate behavior is well approximated by a \textit{single-input single-output, linear time-invariant}  (SISO-LTI) system. Hence the control problem for the balancing authority can be addressed using classical control design methods.  State estimation is not required --- the information required at the BA is an estimate of the proportion of 
loads 
that are operating.  

In the numerical example considered in this paper, the linear system is minimum-phase and stable, which is very helpful for control design.   

The remainder of the paper is organized as follows.   The control solution for a single pool is described in
\Section{s:ppcontrol},  along with approximations of the optimal control solution based on general theory presented in the Appendix.   The control of the aggregate collection of pools is considered in \Section{s:mfg}. Conclusions and directions of future research are contained in \Section{s:conclude}.

\section{Optimal control for a load and for the grid} 
\label{s:ppcontrol}


\begin{figure}[h]
\vspace{-.25cm}
\Ebox{.7}{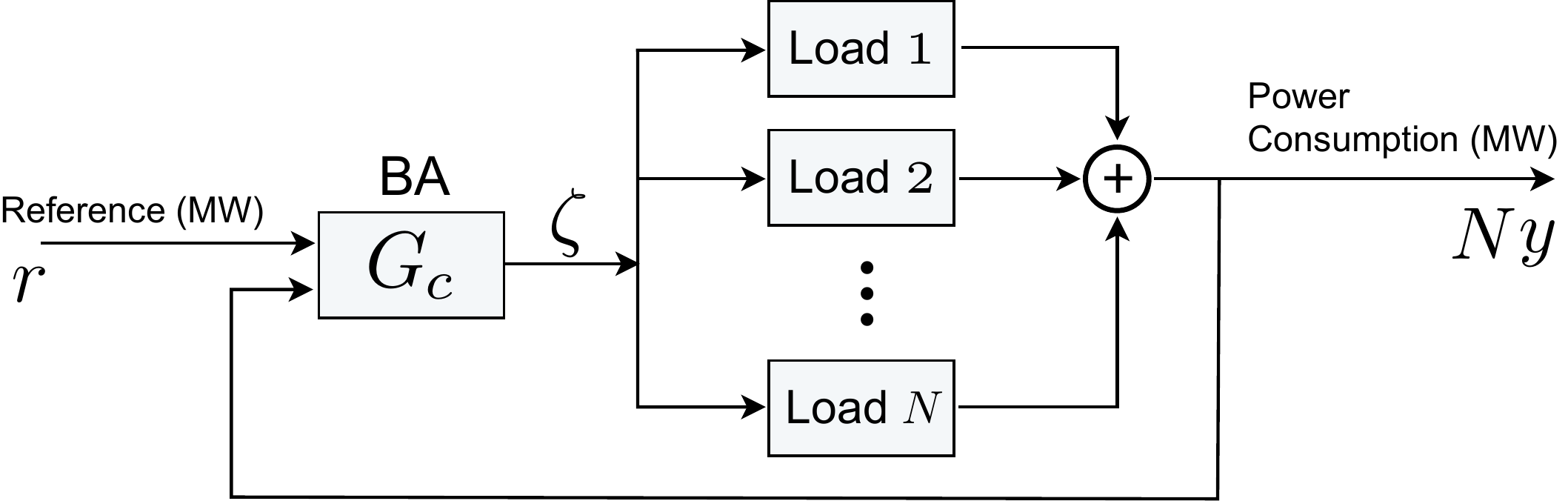}
\vspace{-.25cm}
\caption{The control architecture: command $\bfzeta$ is computed at a BA,
 and   transmitted to each pool pump. The control decision at a load
 is binary (turn on/off), and is based only on its own state and the signal $\bfzeta$.}
\label{fig:arch}
\vspace{-08pt}
\end{figure} 

\subsection{Control architecture overview}

We begin with a description of the control and information architecture that is the subject of this paper.   The components of the architecture are illustrated in  \Fig{fig:arch}:  
\begin{romannum}
\item 
There are $N$ homogeneous loads that receive a common scalar command signal from the balancing authority, or BA, denoted $\bfzeta=\{\zeta_t\}$ in the figure.

Randomization at each load is desirable to avoid synchronization of loads, and also to facilitate analysis of the aggregate system.  It is assumed that
each load evolves as a controlled Markov chain: 
The transition probability for each load is determined by its own state, and  the BA signal $\bfzeta$.  The common dynamics are defined by a controlled transition matrix $\{P_\zeta : \zeta\in\Re\}$.  For the $i$th load, there is a state process $\bfmX^i$ whose transition probability at time $t$ is given by,
\begin{equation}
\Prob\{X^i_{t+1} = x^+ \mid X^i_t = x^- ,\,  \zeta_t=\zeta\} = P_{\zeta}(x^-,x^+) 
\label{e:Pzeta}
\end{equation}
where $x^-$ and $x^+$ are possible   state-values.
The details of the model
are described in \Section{s:loadmodel}.

\item
The BA has measurements of the other two scalar signals shown in the figure:  The normalized aggregate power consumption $\bfmy$ and  desired deviation in power consumption $\bfmr$.

When $\zeta_t=0$ for all $t$, then the aggregate power consumption takes the value $\bfmy^0$.   The goal of the BA is a tracking problem:  Achieve $y_t\approx y^0 +r_t$ for all $t$.   This can be addressed using classical control techniques if the dynamics from $\bfzeta$ to $\bfmtily= \bfmy-\bfmy^0$ can be approximated by an LTI system.

\end{romannum}

The main contributions of this paper are based on the construction of the controlled transition matrix for an individual load,  taking into account potentially conflicting goals:  The BA desires overall dynamics from $\bfzeta$ to $\bfmy$ that facilitate tracking the reference signal $\bfmr$.  Each load requires good  quality of service.  In the case of a pool, the water must be kept clean, and the electricity bill must remain constant over each month.

An approach of Todorov \cite{tod07} is adopted to construct the family of transition matrices $\{P_\zeta : \zeta\in\Re\}$.  They are smooth in the parameter $\zeta$,  and a first-order Taylor series approximation gives, for any pair of states $(x^-,x^+) $
\begin{equation}
P_\zeta(x^-,x^+) = \exp(\zeta \Gamma(x^-,x^+) +O(\zeta^2)) P_0(x^-,x^+) 
\label{e:expclE}
\end{equation}
where $P_0$ denotes the dynamics of a load when $\bfzeta\equiv 0$,  and $\Gamma$ is a matrix.  Based on \eqref{e:clE}, we have
\[
\Gamma(x^-,x^+)
 = \tilutil(x^-)  +\derH (x^+) -\derH (x^-) 
\]
where the function $\derH$ is a solution to \textit{Poisson's equation}; a linear equation for the nominal model.

This structure leads to the   LTI approximation of the input-output dynamics from $\bfzeta$ to $\bfmtily$ that is presented in \Proposition{t:linear}.  \Section{s:approx} also contains second-order approximations of $P_\zeta$.

In \Section{s:mfg} these general techniques are applied to a collection of residential pools.  In this example it is found that the LTI model is minimum phase,  and that a simple PI controller can be effectively used for the control transfer function $G_c$ shown in \Fig{fig:arch}.

\subsection{Load model and design}
\label{s:loadmodel}

In this section we present a procedure to construct the controlled transition matrix appearing in \eqref{e:Pzeta}.  The controlled Markov chain evolves on a finite state space, denoted $\state = \{x^1,\dots,x^d\}$.      The construction is based on an optimal control problem for an individual load, taking into account the needs of the load and the grid.  

It is assumed that a transition matrix $P_0$ is given that models ``control free'' behavior of the Markov chain,  and  a utility function $\util\colon\state\to \Re$ is used to model the needs of the grid.    The optimal control problem will balance average utility and the cost of deviation.  

Since we focus on a single load, in this subsection the index $i$ in  \eqref{e:Pzeta} is dropped,  and we denote by $\bfmX=(X_0,X_1,\dots)$ the stochastic process evolving on $\state$ that models this load.

\begin{figure}[h]
\Ebox{.55}{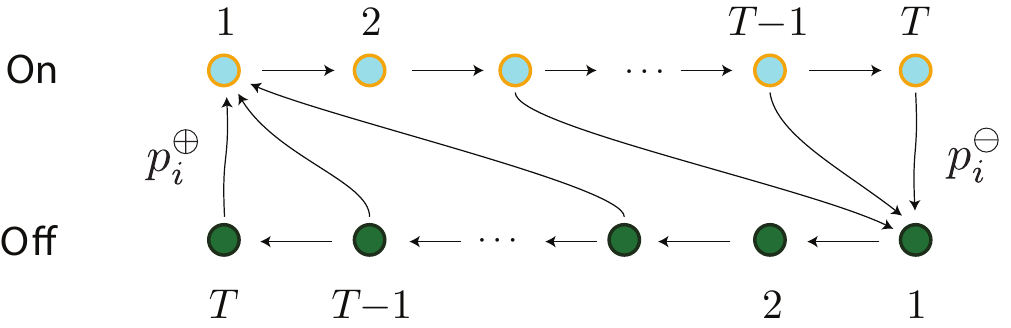} 
\vspace{-.25cm}
\caption{State transition diagram for the pool-pump model.
}
\label{fig:pppDynamics}
\vspace{-08pt}
\end{figure}

In the second half of the paper we will focus on a particular example in which each load is a residential pool pump.
The true nominal behavior would be deterministic -- most consumers set the pump to run a fixed number of hours each day.
However, the randomized policy is based on a stochastic model for nominal behavior, so we introduce some randomness to define the nominal transition matrix
$P_0$.   The state space is taken to be the finite set,
\begin{equation}
\state=\{  (m,i) :  m\in  \{  \oplus,\ominus\}  ,\  i\in \{1,\dots,T\} \}  
\label{e:poolstate}
\end{equation}
If 
$X_t = (\ominus,i)$, this indicates that the pool-pump was turned off and has remained off for $i$ time units, and $X_t = (\oplus,i)$ represents the alternative that the pool-pump has been operating continuously for exactly $i$ time units.  
A state-transition diagram is shown in   \Fig{fig:pppDynamics}.   The values of $P_0(x,y)$ will be chosen to be nearly $0$ or $1$ for most $x,y\in\state$.

The utility function $\util$ on $\state$ is chosen as the indicator function that the pool pump is   operating: 
\[
\util(x) =  \sum_i \ind\{ x = (\oplus, i) \}
\]
Whether this actually represents any utility to the grid operator depends on the state of the grid.  This will be clarified after we define the optimization problem and its solution.

\medbreak

\textit{We now return to the general model.}   Consider first a finite-time-horizon optimization problem:  
For a given terminal time $T$,  let $p_0$ denote the probability measure on strings of length $T$:
\[
p_0(x_1,\dots, x_{T}) =   \prod_{i=0}^{T-1} P_0(x_i,x_{i+1}),\qquad x\in \state^{T}
\]
where $x_0\in\state$ is assumed to be given.

A fixed scalar $\zeta\in\Re$ is interpreted as a \textit{weighting parameter} in the following definition of  \textit{total welfare}.
For any  probability measure $p$,    this is defined as the weighted difference,
\[
\welf_T(p) =  \zeta \Expect_p\Bigl[\sum_{t=1}^{T} \util(X_t) \Bigr] -D(p\| p_0) \,
\]
where the expectation is with respect to $p$,  and $D$ denotes relative entropy.  
Let $p^{T*}$ denote the probability measure that 
maximizes this expression.

\begin{proposition}
\label{t:twisted}
The probability measure $p^{T*}$ is   the \textit{twisted distribution},
\begin{equation}
p^{T*} (x_1,\dots, x_{T}) = \exp\Bigl(\zeta \sum_{t=1}^{T} \util(x_t) - \Lambda_T(\zeta)\Bigr) p_0 (x_1,\dots, x_{T})
\label{e:pstar}
\end{equation}
where 
\begin{equation}
\Lambda_T(\zeta) = \log\Bigl\{ \Expect\Bigl[
\exp\Bigl(\zeta \sum_{t=1}^{T} \util(X_t) \Bigr) \Bigr] \Bigr\}\,,
\label{e:LMGF}
\end{equation}
and the expectation is with respect to $p_0$. 
Moreover, $\welf_T(p^{T*})=\Lambda_T(\zeta) $  is the optimal value.
		\qed
\end{proposition}

\IEEEproof
Optimality of $p^{T*}$   follows from convex duality between the log moment generating function and relative entropy 
-- \cite[Proposition II.1]{huaunnmeyveesur11} and   \cite[Lemma~2.39]{demzei98a}.  The formula \eqref{e:LMGF} follows from the fact that $p^{T*}$ sums to unity, so that $\Lambda_T(\zeta)$ can be interpreted as a normalizing constant.
  
The identity $\welf_T(p^{T*})  = \Lambda_T(\zeta)$ follows from the definitions of $\welf_T$ and $p^{T*}$.
\qed

The probability measure $p^{T*}$ defines a Markov chain on the time
interval $\{0,1, \dots, T\}$, but it is not necessarily time-homogeneous. 
In the
infinite horizon case, we would like to find a distribution $p^*$
on infinite sequences that attains the optimal average welfare,
\notes{SM to PB:  you added the homogeneous statement before the equation, but it already appeared after.}
\begin{equation}
\eta^*_\zeta=\lim_{T\to\infty} \frac{1}{T} \welf_T(p^{T*})  
=
\lim_{T\to\infty} \frac{1}{T}
 \log\Bigl\{ \Expect\Bigl[
\exp\Bigl(\zeta \sum_{t=1}^{T} \util(X_t) \Bigr) \Bigr] \Bigr\}
\label{e:etastar}
\end{equation}
A solution to the infinite horizon problem is  given
by a time-homogenous Markov chain whose transition matrix $\cP_\zeta$ is easy to compute, based on  the solution of an
eigenvector problem;  these results are summarized in the proposition that follows.
The proof of \Proposition{lem:Pcheck-infinite}
is given in the Appendix.

\begin{proposition}
\label{lem:Pcheck-infinite}
If $P_0$ is irreducible, an optimizing $p^*$ that achieves \eqref{e:etastar} is defined by a time-homogeneous Markov chain whose transition probability  is given by
\begin{equation}
\cP_\zeta(x,y) = \frac{1}{\lambda  }  \frac{1}{  v(x)}   \haP_\zeta(x,y) v(y)\,,\qquad x,y\in\state ,
\label{e:cPool}
\end{equation}
where $\haP_\zeta$ is the scaled transition matrix,
\begin{equation}
\haP_\zeta(x,y) = \exp(\zeta \util(x)) P_0(x,y) \,,\qquad x,y\in\state,
\label{e:scaledTM}
\end{equation}
  and $\lambda,v$ is that eigen-pair corresponding to the eigenvector problem
\begin{equation}
\haP_\zeta v = \lambda v
\label{e:TodPF}
\end{equation}
such that $\lambda=\lambda_\zeta>0$ is the unique maximal eigenvalue for $\haP_\zeta$,
$v=v_\zeta$ is unique up to constant multiples,  and $v(x)>0$ for each
$x\in\state$.   In addition, the following bounds hold for each $T$,
\begin{equation}
\begin{aligned}
0\le  \welf_T(p^{T*}) -  \welf_T(\cp^{T}) & \le 2 \| h \|_{\text{\small sp}}
\\
|T\Lambda -  \welf_T(p^{T*})   |& \le  \| h \|_{\text{\small sp}}
\\
\end{aligned}
\label{e:cpOpt}
\end{equation}
where $h=\log v$, and the span norm is defined by $ \| h \|_{\text{\small sp}} = \max h - \min h$.  

Consequently,  the Markov model achieves the optimal average welfare \eqref{e:etastar} with $\eta^*_\zeta=\Lambda$. 
		\qed
\end{proposition}

\notes{
 $\Lambda=\log(\lambda)$PB: is this
  equality a definition or a result? SM: result - clearer now?} 

The  eigenvector problem~\ref{e:TodPF} appears in multiplicative ergodic theory \cite{konmey05a}, and also in
Todorov's analysis \cite{tod07}. It is shown in  \cite{tod07} that 
the \textit{relative value function} appearing in the average cost optimality equations is the logarithm of the eigenvector:
\begin{equation}
h^*(x) = \log(v(x)),\qquad x\in\state.
\label{e:relative}
\end{equation}
See also the derivation in  \cite{meybarbusehr13} for a variant of this model.

Second order Taylor series approximations for $v$ and $\eta^*$ near $\zeta\approx 0$ can be found by borrowing tools from large-deviations theory.  Some of these approximation results are new, and are collected together in the next section and in the Appendix.

\subsection{Approximations}
\label{s:approx}

Approximations will be needed for analysis when we extend the model to allow $\zeta$ to change with time.

A solution to the eigenvector problem \eqref{e:TodPF} can be represented through a regenerative formula.   Let $\atom\in\state$ be some fixed state that is reachable from each initial condition of the chain, under the transition law $P_0$.  That is, the chain is assumed to be $\atom$-irreducible \cite{MT}.  Since the state space is assumed to be finite, it follows that there is a unique invariant probability measure $\pi_0$ for $P_0$.    The first  return time is denoted,
\[ 
	 \tau = \min\{ t\ge 1 : X_t = \atom\}\,.
\]

Recall that the infinite horizon optimal welfare  is given by $\eta^*_\zeta=\log(\lambda)$.
From the theory of positive matrices \cite{sen81,num84,konmey05a}, it follows that it is the unique solution to,
\begin{equation}
1=
\Expect_\atom \Bigl[\exp \Bigl(\sum_0^{\tau-1} [ \zeta \util(X_t) - \eta^*_\zeta ] \Bigr)\Bigr]
\label{e:PFeta}
\end{equation} 
where the subscript indicates that the
   initial condition is $X(0)=\atom$.
Moreover, for each $x\in\state$,  the value of $v(x)$ is obtained as the expected sum,
with initial condition $X(0)=x$:
\begin{equation}
v(x)=  \Expect_x\Bigl[\exp \Bigl(\sum_0^{\tau -1} [ \zeta \util(X_t) - \eta^*_\zeta ] \Bigr)  \Bigr]
\label{e:PFv}
\end{equation}
These expectations are each with respect to the nominal transition law $P_0$.

A Taylor-series approximation of $\eta^*_\zeta$
is based on two parameters, defined with respect to the nominal model $P_0$ with invariant probability measure $\pi_0$.  
The first-order coefficient is the
the  steady-state mean of $\util$,
\begin{equation}
\eta_0=\sum_x \pi_0(x)\util(x)
\label{e:eta0}
\end{equation}
The second-order coefficient is based on the  \textit{asymptotic variance} of $\util$ for the nominal model (the variance appearing in the Central Limit Theorem (CLT) for the nominal model). For this finite state space model this has two similar representations,
\begin{equation}
\begin{aligned}
\avar
	&=\lim_{T\to\infty}\frac{1}{T} \Expect\Bigl[ \Bigl( \sum_0^{T -1}    \tilutil(X_t) \Bigr)^2\Bigr]
\\
	&= \pi_0(\atom) \Expect_\atom\Bigl[ \Bigl( \sum_0^{\tau -1}    \tilutil(X_t) \Bigr)^2\Bigr]
\end{aligned}
\label{e:varsigma}
\end{equation}
where $\tilutil = \util-\eta_0$.  See  \cite[Theorem~17.0.1]{MT} for the CLT,  and eqn.~(17.13) of \cite{MT} for the second representation above.

Similarly, the following functions of $x$   are used to define a second order Taylor series approximation for $ h^*_\zeta$.
The first-order term is the solution to \textit{Poisson's equation} for $P_0$,
\begin{equation}
\derH (x) =  \Expect_x\Bigl[\sum_0^{\tau -1}    \tilutil(X_t)   \Bigr]
\label{e:fish0}
\end{equation} 
The asymptotic variance  can be expressed in terms of Poisson's equation \cite{MT,CTCN}:
\[
\avar =\sum_x  \pi_0(x) \bigl(2\tilutil(x)\derH (x) -\tilutil(x)^2 \bigr)
\] 
The second-order term in an approximation of $v$ is another variance,
\begin{equation}
\varH(x)=  \Expect_x\Bigl[  \Bigl(\sum_0^{\tau -1} \tilutil(X_t)  \Bigr)^2 \Bigr] -  \bigl(\derH (x)\bigr)^2 \,,  \quad x\in\state.
\label{e:SOvH}
\end{equation}

\begin{proposition}
\label{t:etaSecondOrder}
The following hold for the finite state space model in which $P_0$ is irreducible:
\begin{romannum}
\item 
The optimal average welfare $\eta^*_\zeta$ is convex as a function of $ \zeta$, and admits the Taylor series expansion, 
\begin{equation}
\eta^*_\zeta
=
\eta_0 \zeta + \half \avar  \zeta^2 + O( \zeta^3)
\label{e:SOeta}
\end{equation}

\item
The mean of $\util$ under the the invariant probability measure  $\cpi_\zeta$ for $\cP_\zeta$ is given by,
\begin{equation}
\sum_x \cpi_\zeta(x)\util(x)
= 
	\frac{d}{d\zeta} \eta^*_\zeta 
\label{e:ppoffProb}
\end{equation}
This admits the first-order Taylor series approximation
\begin{equation}
\frac{d}{d\zeta} \eta^*_\zeta 
=
 \eta_0 +   \avar \zeta + O( \zeta^2)
\label{e:ppoffProbApprox}
\end{equation}

\item
The relative value function \eqref{e:relative}
 admits the second-order Taylor series approximation,
\begin{equation}
 h^*_\zeta(x)
=
 \zeta \derH (x) + \half   \zeta^2 \varH(x) + O( \zeta^3)
\label{e:happrox}
\end{equation}
\end{romannum} 
		\qed
\end{proposition}

\IEEEproof
Equations \eqref{e:SOeta}---\eqref{e:ppoffProbApprox} follow from the fact that $\eta^*_\zeta = \log(\lambda)$ can be expressed as a cumulative log-moment generating function \cite[Prop. 4.9]{konmey05a}.

Convexity follows from the fact that $\eta^*_\zeta$ is the maximum of linear functions of $ \zeta$ (following the linear-program formulation of the ACOE \cite{bor02a}).  
 
 \notes{The following also might be obtained from extending \cite[Prop. 4.9]{konmey05a}.  I don't have a reference, but it is obvious from the representation of $v$}
 
The approximation  \eqref{e:happrox} follows from the  the representation \eqref{e:PFv} for $v$,
and the definition  $h^*_\zeta =  \log(v)$ (see \eqref{e:relative}).
\qed

The representations in this subsection are useful for analysis, but not for computation.   Methods to compute  $\derH $ and $\varH$ are contained in the Appendix. 

\subsection{Aggregate load model}

Consider $N$ loads operating independently under   the randomized policy described in the previous section.  The state of the $i$th load is denoted $ X^i_t  $.
For large $N$ we have from the Law of Large Numbers,
\begin{equation}
\begin{aligned}
\frac{1}{N}\sum_{i=1}^N  \util(X^i_t) & \approx   \Expect [\util (X_t ) ] 
\end{aligned}
\label{e:AvgCostPool}
\end{equation} 
The expectation and probability on the right are with respect to the optimal transition law $\cP_\zeta$, where $ \zeta$ is the parameter used in \eqref{e:etastar}.

We pose the following centralized control problem:  How to choose the
variable $ \zeta$ to regulate average utility \textit{in real time},
based on measurements of the average utility,  and also a regulation
signal denoted $\bfmr$.  Let $y_t$ be the fraction of loads that are
on at time $t$:
\begin{equation}
y_t =\frac{1}{N}\sum_{i=1}^N  \util(X^i_t),
\label{e:y}
\end{equation}
which is assumed to be observed by the BA.

To address the control problem faced by the BA it is necessary to relax the assumption that this parameter is fixed.  
We let $\bfzeta=\{\zeta_0,\zeta_1,\dots\}$ denote a sequence of scalars, which  is regarded as an input signal for the control problem faced by the BA.  An aggregate model is obtained in two steps.

In step 1   the existence of
a mean-field limit is assumed: 
Let $N\to\infty$ to obtain the generalization of \eqref{e:AvgCostPool},
\begin{equation}
\lim_{N\to\infty}
\frac{1}{N}\sum_{i=1}^N  \ind\{ X^i_t =x \} = \mu_t(x)\,, \quad x\in\state.  
\label{e:mfgPool}
\end{equation}  
For a given initial distribution $\mu_0$ on $\state$, 
the distribution $\mu_t$   is defined by  $\mu_t(x_t) = {}$
\begin{equation}
 \sum_{x_i\in\state} \mu_0(x_0) \cP_{\zeta_0}(x_0,x_1) \cP_{\zeta_1}(x_1,x_2) \cdots  \cP_{\zeta_{t-1}}(x_{t-1},x_t)
\label{e:muMF}
\end{equation}
where $x_t$ is an arbitrary state in $\state$, and the sum is over all intermediate states. 
We view $\{\mu_t\}$ as a state process that is under our control through $\bfzeta$.  

Justification for the mean-field limit is contained in
\Theorem{t:MFL}.

Step 2 is based on the Taylor series approximations  surveyed in the previous section to  approximate this nonlinear system by a linear state space model with $d$-dimensional state $\bfPhi$ and output $\bfgamma$. 
It is defined so that for any time $t$, and any $i$,
\[
\begin{aligned}
\mu_t(x^i)&=\pi_0(x^i)+\Phi_t(i)   + o(\bfzeta)
\\
\gamma_t  &= \tily_t + o(\bfzeta)
\end{aligned}
\]
where 
$ \tily_t=y_t-y^0$, with $y^0 = \sum_x \pi_0(x) \util(x)$, and
 where $o(\bfzeta)$ is in fact $O(\zeta_0^2+\cdots + \zeta_t^2)$.   \notes{PB: why is the IC part of the
  linearization description? SM:  ok now?  Remember, some readers might be outside of control area. }

\begin{proposition}
\label{t:linear}
Consider the nonlinear state space model whose state evolution is
$\mu_{t+1}
= \mu_t \cP_{\zeta_t}$,  and output is $y_t=\sum_x \mu_t(x)\util(x)$.  
Its unique equilibrium with $\bfzeta\equiv 0$ is $\mu_t\equiv \pi_0$ and $y_t\equiv y^0\eqdef \sum_x \pi_0(x)\util(x)$.   Its linearization around this equilibrium is given by,
\begin{equation}
\begin{aligned}
 \Phi_{t+1} &=  A \Phi_t + B \zeta_t
 \\
\gamma_t &=  C \Phi_t
\end{aligned}
\label{e:LSSmfg}
\end{equation}
where
$A=P^\transpose_0$,   $C$ is a row vector of dimension $d=|\state|$ with $C_i= \util(x^i) $ for each $i$, and $B$ is a $d$-dimensional column vector with entries $B_j = \sum_x\pi_0(x) \clE(x,x^j) $,  where
\begin{equation}
\clE(x^i,x^j)
 = \Bigl[  \tilutil(x^i)+ \derH (x^j) -\derH (x^i)  \Bigr]P_0(x^i,x^j) 
\label{e:clE}
\end{equation}
for each $x^i,x^j\in \state$. The initial condition is
$\Phi_0(i)=\mu_0(x^i)-\pi_0(x^i)$,  $1\le i\le d$.		

The matrix $\clE$ is equal to the derivative,
\[
\clE=\frac{d}{d\zeta} P_\zeta \Big|_{\zeta=0}
\]
Consequently,  the formula \eqref{e:clE} implies the approximation \eqref{e:expclE}.
\qed
\end{proposition}
%

\IEEEproof
The formulae for $A$ and $C$ follow from the fact that the system is linear in the state.  We have, from \eqref{e:cPool},  
\[
\cP_\zeta(x^i,x^j) = e^{ \zeta\util(x^i) - \eta^*_\zeta - h^*_\zeta(x^i)} P_0(x^i,x^j) e^{h^*_\zeta(x^j)}
\]
Based on the first order approximation of $h^*_\zeta$ in  \Proposition{t:etaSecondOrder}  we obtain, 
\[
\cP_\zeta(x^i,x^j) \approx e^{\zeta[-\derH (x^i)+\tilutil(x^i) ]} P_0(x^i,x^j) e^{\zeta \derH (x^j)}
\]
where $\derH $ is a solution to Poisson's equation (with forcing function $\util$) for the nominal model  (see \eqref{e:fish0}).
Using a first order Taylor series for the exponential then gives,
\[
\begin{aligned}
\cP_\zeta(x^i,x^j) 
&\approx [1-\zeta(\derH (x^i)-\tilutil(x^i) )] P_0(x^i,x^j)[1+ \zeta \derH (x^j)] 
\\
&\approx P_0(x^i,x^j) + \zeta \clE(x^i,x^j)
\end{aligned}
\]

If $\mu \approx \pi_0$
and $ \zeta$ is small, then we can approximate,
\[
\mu \cP_\zeta \approx \mu P_0 + \zeta B^\transpose  \,,
\]
where $B$ is the column vector with entries $B_j = \sum_x\pi_0(x) \clE(x,x^j) $.
\qed

Next we justify the mean-field model \eqref{e:mfgPool}.

For the purpose of analysis we lift the state space from the $d$-element set
$\state = \{x^1,\cdots,x^d\}$,  to the $d$-dimensional simplex $\Spx$. 
For the $i^{th}$ load at time $t$, the element $\pi_t^i \in \Spx$ is the degenerate distribution whose mass is concentrated at $x$ if   $X^i_t= x$. The average over $N$, denoted $\mu_t^N\in \Spx$,  is the empirical distribution,
\[	 
\mu_t^N( x)
=\frac{1}{N}\sum_{i=1}^N  \pi_t^i(x)  
=\frac{1}{N}\sum_{i=1}^N  \ind\{X^i_t =x \} 
\, , \quad x\in\state,
\]
In the proof of convergence it is assumed that $\bfzeta^N$ is obtained using state feedback of the form,
\[
\zeta_t^N = \phi_t(\mu_0^N,\dots,\mu_t^N)
\]
where $\phi_t\colon\Spx^{t+1}\to\Re$ is continuous for each $t$, and
does not depend upon $N$. The following result establishes convergence.

\begin{theorem}
\label{t:MFL}
Suppose $\mu_0^N \rightarrow \mu_0$  as  $ N\rightarrow \infty$, and that the state transition matrix $P_\zeta$ is   continuous as a function of $\zeta$.  Then for each $t$,
\begin{equation}
	\lim_{N\to\infty}  \mu_t^N= \mu_t, \qquad \text{\it with probability one, }
\label{e:MFL}
\end{equation}
where the right hand side denotes the probability measure \eqref{e:muMF},  in which
\[
\zeta_t = \phi_t(\mu_0,\dots,\mu_t),\qquad t\ge 0.
\]
\qed
\end{theorem}

The proof of this result is given at the end of this subsection, and is largely based on a version of the Law of Large Numbers. 

Let $\{M_{N,k},1\le k \le N\}$ denote a martingale array:  
This means that $\Expect[M_{N,k}|M_{N,1},\cdots M_{N,j}]=M_{N,j}$ for each $N$ and $1 \le j < k \le N$. When $k=N$, we denote $M_N = M_{N,N}$.
\begin{proposition}
\label{t:MA-LLN}
Suppose that $M_{N,k}$ is a martingale array with bounded increments:  For some $c_m<\infty$,
\[
 	| M_{N,k+1}-M_{N,k} | \le c_m\qquad \text{for all $k$ and $N$}
\]
Then the Law of Large Numbers holds: 
\[
	\lim_{N\to\infty} \frac{M_N}{N} = 0, \qquad \text{\it with probability one. }
\]
\qed
\end{proposition} 

\IEEEproof
The Hoeffding-Azuma inequality
\cite{mcd98a} gives the following bound:
\[
\Prob\{ N^{-1} |M_N |\ge t\}  \le 2 \exp(- [N t]^2/[2 N c_m^2] )
\]
The right hand is summable, so the result follows from the Borel-Cantelli Lemma.
\notes{Isn't this an amazing bound!?}
\qed 

\Proposition{t:MA-LLN} is applied to show that the sequence of empirical distribution $\mu_t^N$ can be approximated by the mean-field model perturbed by a  disturbance that vanishes as $N\to\infty$:

\begin{lemma}
\label{t:W-is-MA}
The empirical distributions $\{\mu_t^N: t\ge 0\}$ obey the recursion  
\begin{equation}
	 \mu_{t+1}^N = \mu_t^N P_{\zeta_t^N}+W_{t+1}^N,
	\label{e:empir_dist}
\end{equation}
in which,  $W_{t+1}^N=\frac{1}{N}\sum_{i=1}^N \Delta_{t+1}^i$ for a family of vector random variables $\{ \Delta_{t+1}^i\}$.

On denoting $M_{N,k}=\sum_{i=1}^k \Delta_t^i$ we have,
\begin{romannum} 
\item $\{M_{N,k}:1\le k\le N\}$ is a martingale array.
\item There exits $c_m$ such that $\| M_{N,k}-M_{N,k-1}\| \le c_m$  for all $N$ and all $k$ such that $1< k\leq N$.
\end{romannum}
\end{lemma}

\textit{Proof of \Lemma{t:W-is-MA}}: 
To establish \eqref{e:empir_dist} we first establish a similar expression for $\{\pi_t^i \}$.

For each $i$, the sequence of degenerate  distributions $\{\pi_t^i \}$ evolve according to a random linear system,
\begin{equation}
	\pi_{t+1}^i=\pi_t^iG_{t+1}^i 
\label{e:piG}
\end{equation}
in which $\pi_t^i$ is interpreted as a $d$-dimensional row vector, and $G_t^i$ is a $d\times d$ matrix with entries $0$ or $1$ only, and $\sum_l G_t^i(x^j,x^l)=1$ for all $j$. It is conditionally independent of $\{\pi_0^i,\cdots,\pi_t^i\}$, given $\zeta^N_t$, with
\begin{equation}
\label{e:EG=P} 
 	\Expect[G_{t+1}^i|\pi_0^i, \cdots, \pi_t^i, \zeta^N_t]=P_{\zeta^N_t}.
\end{equation} 
Dependency of $\pi_t^i$, $G_t^i$ on $N$ is suppressed, but we must distinguish $\zeta^N_t$ from its limit $\zeta_t$.

The random linear system \eqref{e:piG} can thus be described as a linear system driven by ``white noise'':
\begin{equation}
	\pi_{t+1}^i=\pi_t^iP_{\zeta^N_t}+\Delta_{t+1}^i
	\label{e:Sys_delta}
\end{equation}
where, $\{\Delta_{t+1}^i=\pi_t^i(G_{t+1}^i -P_{\zeta^N_t}): t\geq1 \}$, 
which establishes  \eqref{e:empir_dist}.

The following representation will clarify the remaining analysis:
\notes{SM:  I need to make its use more apparent}
\begin{equation}
\text{ $G_t^i=\clG(\zeta^N_{t-1},\xi_t^i)$, where $\{ \xi_t^i: t\geq 1, \ i\ge 1\}$ are i.i.d.}
\label{e:A1}
\end{equation}
For $1\le i< N$ and fixed $t$, we define  two $\sigma$-algebras:
\[
\begin{aligned}
                \clF_i &= \sigma \{\Delta _t^k, k\le i \}
                \\
               \clH_i&=\sigma\{\pi_{t-1}^{k+1}, \zeta^N_{t-1}, \Delta _t^k, k\le i  \}
\end{aligned}
\] 
Under \eqref{e:A1} we have the   extension of $\eqref{e:EG=P}$, that $\Expect[G_t^{i+1}\mid \clH_i]=P_{\zeta^N_{t-1}}$. 
Moreover, by construction the random variable $\pi_{t-1}^{i+1}$ is $\clH_i$-measurable. 
Therefore,
\[
	\Expect[\Delta_t^{i+1}\mid \clH_i]=\Expect[\pi_{t-1}^{i+1}(G_t^{i+1}-P_{\zeta^N_{t-1}})\mid \clH_i]=0
\]
The smoothing property of the conditional expectation, and the construction $ \clF_i \subset \clH_i$, then gives (i),
\[
		\Expect[\Delta_t^{i+1}\mid \clF_i]=\Expect[\Expect[\Delta_t^{i+1}\mid \clH_i]\mid \clF_i]=0
\]

From the definition of $\Delta_t^i$ below equation \eqref{e:Sys_delta}, it follows that $\{ \|\Delta_t^i\| \}$ admits a uniform bound. Consequently, $\| M_{N,k}-M_{N,k-1}\|=\| \Delta_t^k\|$ is bounded, which is (ii).
\qed

\head{Proof of \Theorem{t:MFL}}

Denote, for $T\ge 0$, the deviation  $ \tilde{\mu}_T^N   = \mu_T^N  -\mu_T$.  

We prove   by induction on $T$ that $\tilde{\mu}_T^N \to 0$ as $N\to\infty$.
This holds by assumption when $T=0$.

Suppose now that \eqref{e:MFL} holds for $t\le T$. 
By continuity of $\phi_t$, it follows that $\zeta_t^N\to \zeta_t$ as $N\to\infty$.
We also have by the definitions,
\[
 \tilde{\mu}_{T+1}^N 
 	=
	\tilde{\mu}_T^NP_{\zeta_T}
	+ \mu_T^N(P_{\zeta_T^N}-P_{\zeta_T})	+	W_{T+1}^N 
\]

\Lemma{t:W-is-MA} and \Proposition{t:MA-LLN} imply that $W_{T+1}^N 
\to0$ as $N\to\infty$.  Continuity of $P_\zeta$ then implies that 
\[
 \lim_{N\to\infty} \tilde{\mu}_{T+1}^N = 0
 \]
\qed

\section{Controlling a large number of pools} 
\label{s:mfg}

For the remainder of the paper we apply the results of the previous section to the control of a large population of residential pools.   The nominal transition matrix $P_0$  is defined by the probabilities
of turning the pump on or off, as illustrated in the state transition
diagram \Fig{fig:pppDynamics}.   In many of the numerical results described
below a symmetric model was chosen for $P_0$ in which
$p_i^\oplus=p_i^\ominus$, where $p_i^\oplus \eqdef \Prob\{\text{pump
  switches 
  on} \,|\, \text{it has been off  $i$ hours} \}$. Similarly, $p_i^\ominus \eqdef \Prob\{\text{pump
  switches 
off}  \,|\,  \text{it has been on   $i$ hours} \}$.

The utility function $\util$ on $\state$ is chosen as the indicator function that the pool pump is   operating: 
\begin{equation}
\util(x) =  \sum_i \ind\{ x = (\oplus, i) \}
\label{e:kappaOff}
\end{equation}
The parameter $ \zeta$ in \eqref{e:etastar} can be positive or negative;  
If $ \zeta>0$ this control formulation is designed to provide incentive to turn pumps on.

\begin{figure}[h]
\Ebox{.5}{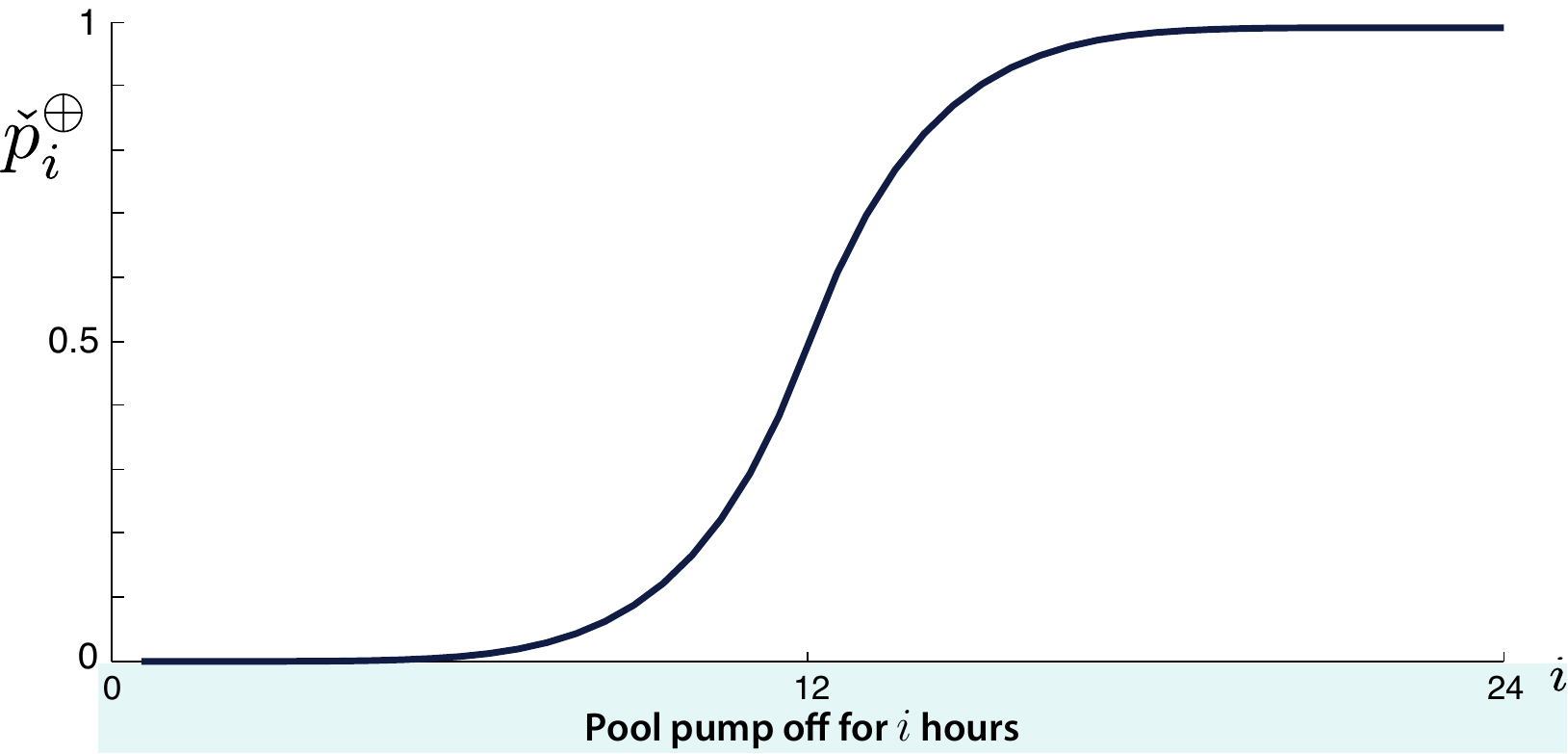}
\vspace{-.05cm}
\caption{Control free behavior of a pool used for numerical studies.}
\label{fig:plus}
\vspace{-08pt}
\end{figure}

It remains to give numerical values for $p_i^\oplus$ and $p_i^\ominus$,  $1\le i\le T$.
In the symmetric model, the specification of these probabilities is performed as follows.
Fix $\gamma>1$ and define,
\[
\varrho_s(x) = \begin{cases}
2^{\gamma-1} x^\gamma &  0\le x\le 1/2
\\
1- 2^{\gamma-1}(1- x)^\gamma &  1/2\le x\le 1
\end{cases}
\]
If over a 24 hour day we choose a sampling time $T=30$ minutes, then in the symmetric model we take,
\begin{equation}
p_i^\oplus = p_i^\ominus =\varrho_s(i/48)\,,\qquad 1\le i\le 48.
\label{e:pvarrho}
\end{equation}
\Fig{fig:plus} shows a plot of the resulting probability $p_i^\oplus$ vs.\ $i$ with $\gamma = 6$.  

To go beyond the asymmetric model,  introduce a parameter $\alpha $ intended to represent the fraction of the day that the pool is operating.  We modify  $ \varrho_s$ as follows,
\[
   \varrho_s^+(x) = \varrho_s(x^{\delta_+}),\qquad  \varrho_s^-(x) =\varrho_s(x^{\delta_-}), 
\]
where $\delta_+$ is chosen so that $0.5^{\delta_+}=1-\alpha$,  or $\delta_+ =- \log_2(1-\alpha)$,
\spm{Please check my work}
and similarly $\delta_- =- \log_2(\alpha)$. For the same sampling parameters as in the previous example, we then take,
\begin{equation}
p_i^\oplus =\varrho_s^+(i/48),\quad p_i^\ominus =\varrho_s^-(i/48)\,,\qquad 1\le i\le 48.
\label{e:pdefn}
\end{equation}
As $\gamma\to\infty$,  the functions in \eqref{e:pdefn} will converge to step functions corresponding to a deterministic cleaning period of $\alpha \times 24$ hours.    
We find numerically that the average cleaning period is somewhat smaller when $\alpha <\half$ and   $\gamma<\infty$.

\subsection{Approximations}

\begin{figure}[h]
\Ebox{.5}{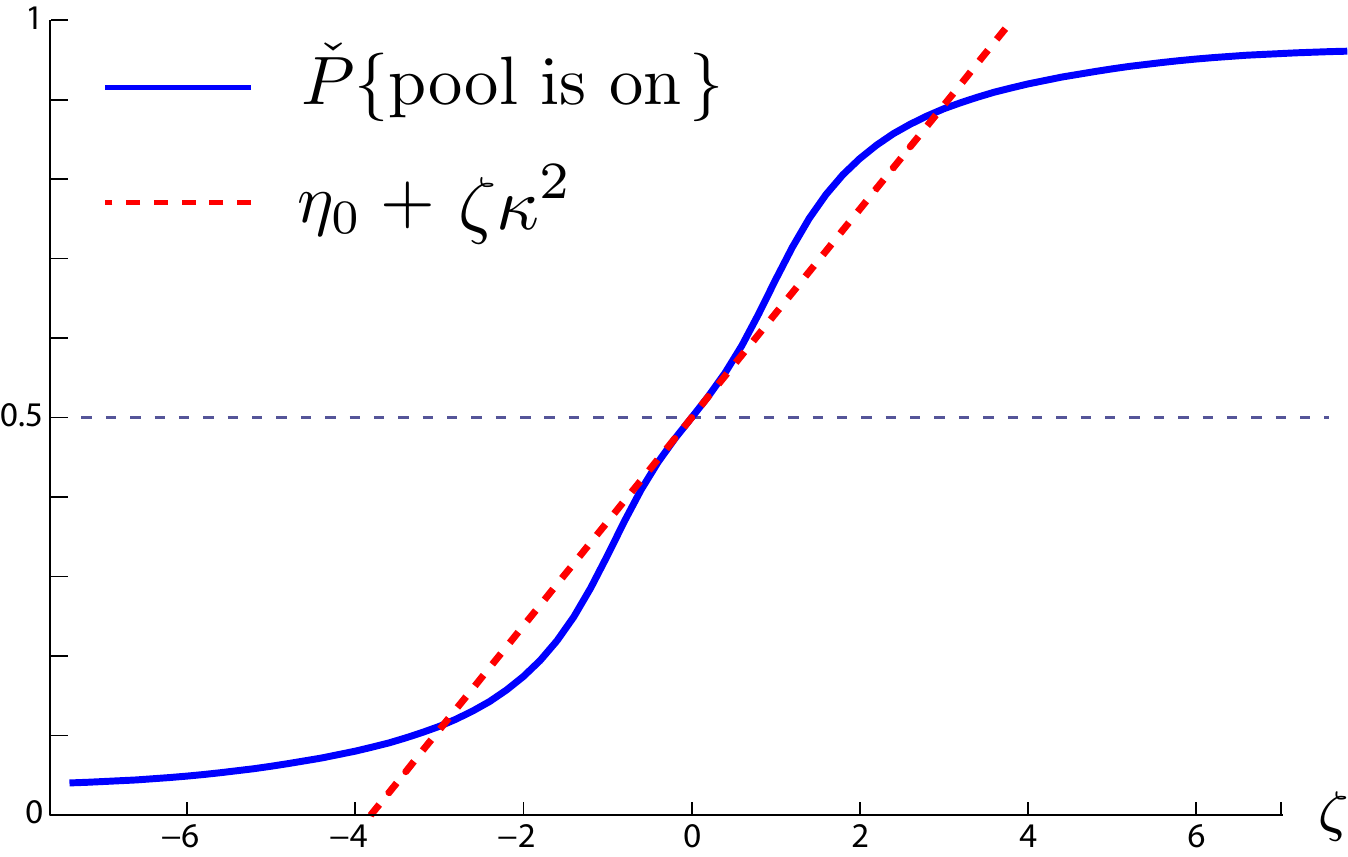}
\vspace{-.2cm}
\caption{Approximation of the steady-state probability that a pool-pump is operating under $\cP$.}
\label{fig:pponprob} 
\end{figure}

The steady-state probability that a pool-pump is in operation is 
given by 
\[
\check P\{\text{pool-pump is on}\}=
 \sum_x \cpi_\zeta(x)\util(x) 
\]
A linear approximation is obtained in \Proposition{t:etaSecondOrder}~(ii):  
\begin{equation}
\check P\{\text{pool-pump is on}\} =   \eta_0 +  \avar \zeta + O( \zeta^2)
\label{e:ppoffProbApp}
\end{equation}     
A comparison of the true probability and its affine approximation   is shown in
\Fig{fig:pponprob} for the symmetric model, in which $\eta_0=1/2$. The approximation is very tight for $|z|\le 3$.  For larger values of $ \zeta$ the true steady-state probability saturates (approximately $0.9$ as $ \zeta\to +\infty$).  

\notes{Yue,  are all the numerics obtained with $\gamma=6$?
\\ Note that it is impossible to approach the values  $\check P\{\text{pool-pump is on}\} =1$ or $0$ with this model.
\\
conclusion:  stick to $|z|\le 3$}

\begin{figure}[h]
\Ebox{.65}{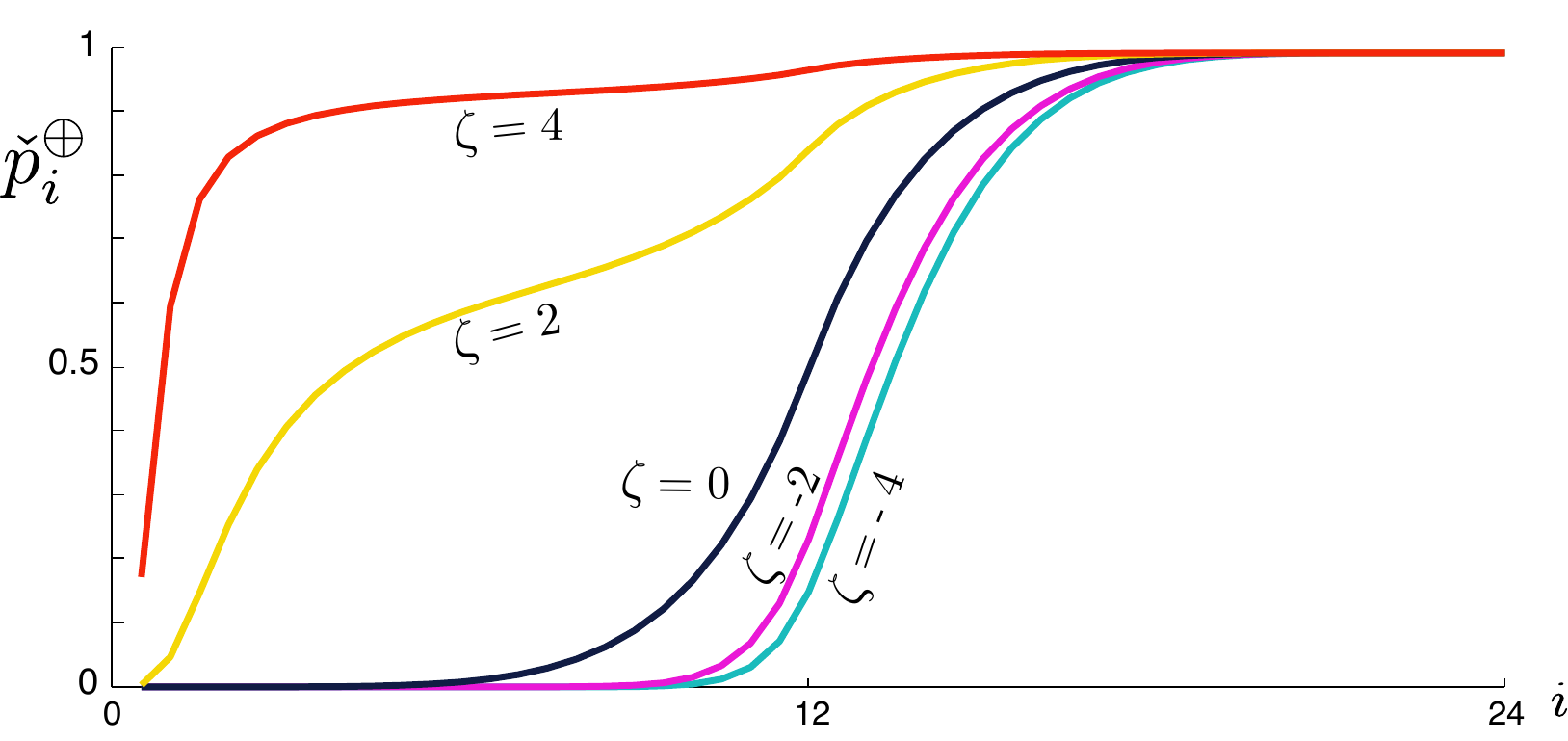}
\vspace{-.2cm}
\caption{Transformed probability vector $\cp^\oplus$ under $\cP$.}
\vspace{-.2cm} 
\label{fig:checkpplus}
\end{figure}

For fixed $ \zeta$,  the controlled model $\cP$ has the same form as $P_0$, with  transformed probability vectors $\cp^\oplus$ and $\cp^\ominus$.   \Fig{fig:checkpplus} contains plots of the transformed vector $\cp^\oplus$ for values $ \zeta=0, \pm 2, \pm 4$.
The plots of  $\cp^\ominus$ are obtained through symmetry.

\begin{figure}[h]
\Ebox{.55}{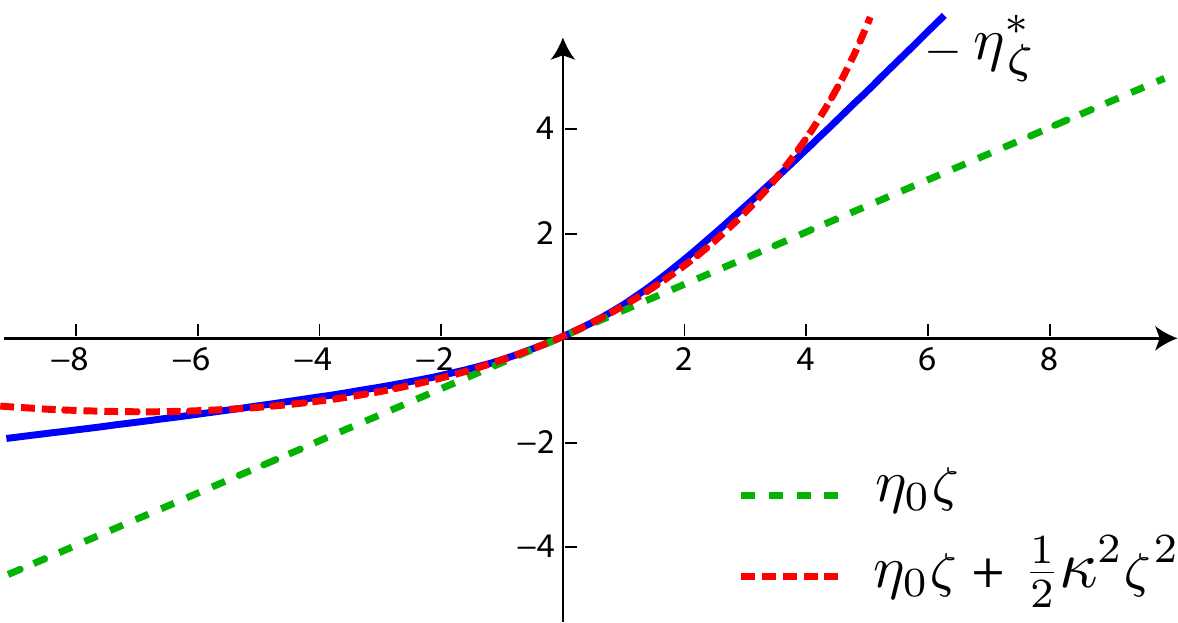}
\vspace{-.25cm}
\caption{The optimal average welfare $\eta_\zeta^*$ and its quadratic approximation.}
\label{fig:QuadApprox} 
\end{figure}

The approximation of the average welfare established in
\Proposition{t:etaSecondOrder}   is,
\begin{equation}
\eta^*_\zeta
=
\eta_0 \zeta+\half \avar  \zeta^2 +O( \zeta^3)
\label{e:SOeta_pool}
\end{equation}
Shown in \Fig{fig:QuadApprox} is a comparison of $\eta_\zeta^*$ with linear and  quadratic approximations based on \eqref{e:SOeta_pool}.

\begin{figure}[h]

\Ebox{.75}{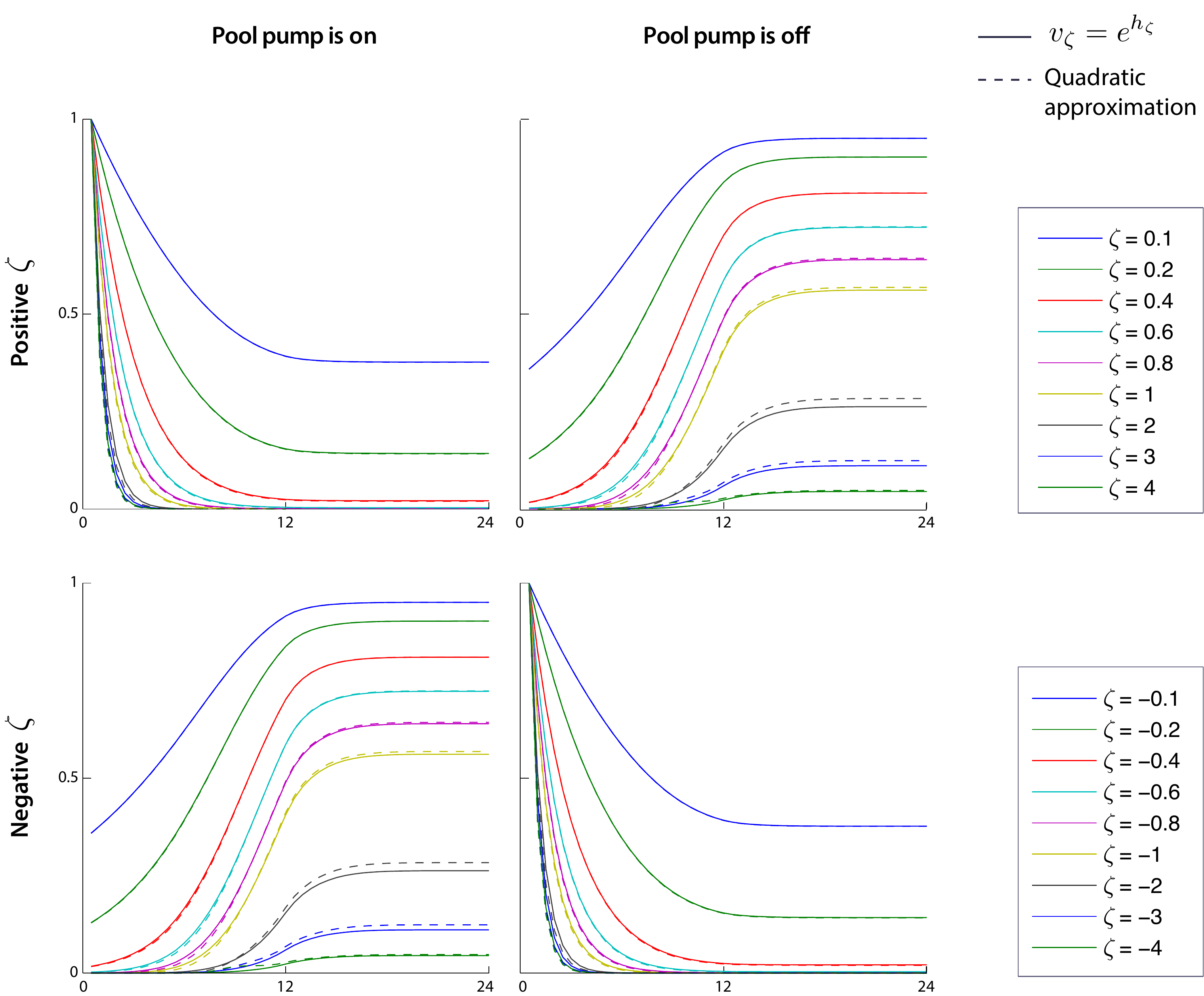} 
\caption{Eigenvectors $v_\zeta=e^{h^*_\zeta}$,
 and their quadratic approximations
$\exp( \zeta \derH (x) + \half   \zeta^2 \varH(x)) $.  }
\label{fig:QuadUapproxExp}
\vspace{-13pt}
\end{figure}

The plots in \Fig{fig:QuadUapproxExp} compare the eigenvector $v=e^{h^*_\zeta}$ with the exponential of the quadratic approximation
\eqref{e:happrox}
given in  \Proposition{t:etaSecondOrder}~(iii).  The computations of $\derH$ and $\varH$ were based on the alternate expression for these functions that are described in \Proposition{t:PoissonDerivatives}.  They are normalized so that the common maxima are equal to unity.   The approximation is nearly perfect for the range of $\zeta\in [-4,4]$.

\subsection{Aggregate load model for pool population}
  

Here we examine the linear model 
\eqref{e:LSSmfg} that will be used by the BA for control synthesis.

We begin with an equilibrium analysis in which $\bfzeta$ is held constant:  
Suppose that $\bfzeta$ does not vary with time, $\zeta_t= \zeta^*$ for all $t$, and consider the steady-state behavior of the  mean-field model.   We denote $y_\infty = \lim_{t\to\infty} y_t$, which is the steady-state probability that a pool is on, for the model with transition law $\cP_{ \zeta^*}$.  This can be approximated using \Proposition{t:etaSecondOrder}:
\[
y_\infty = \Prob\{ \text{Pump is operating}  \}   \approx   \eta_0 +  \avar  \zeta^* 
\] 

From the viewpoint of the BA,  there is a value $G^*$ of desired consumption by all the pools.   If  $\GNoplus>0$ denotes the consumption of one pool pump in operation,  and if there are $N$ pools in total, then the desired steady-state probability is   
$y_\infty = G^*/(N \GNoplus )$.  This translates to a corresponding value of $ \zeta^*$,
\begin{equation}
 \zeta^* 
 \approx  
 \frac{1}{\avar} \Bigl[ \frac{1}{\GNoplus} \frac{G^*}{N} - \eta_0 \Bigr]  
 =
 \frac{1}{\avar}  \frac{1}{\GNoplus} \frac{\widetilde G}{N}  
\label{e:zGapprox}
\end{equation}
where $\widetilde G=G^*-G_0$, with $G_0 = \GNoplus N  \eta_0 $,   the control-free value obtained with $ \zeta^*=0$.

\begin{figure}
\Ebox{1}{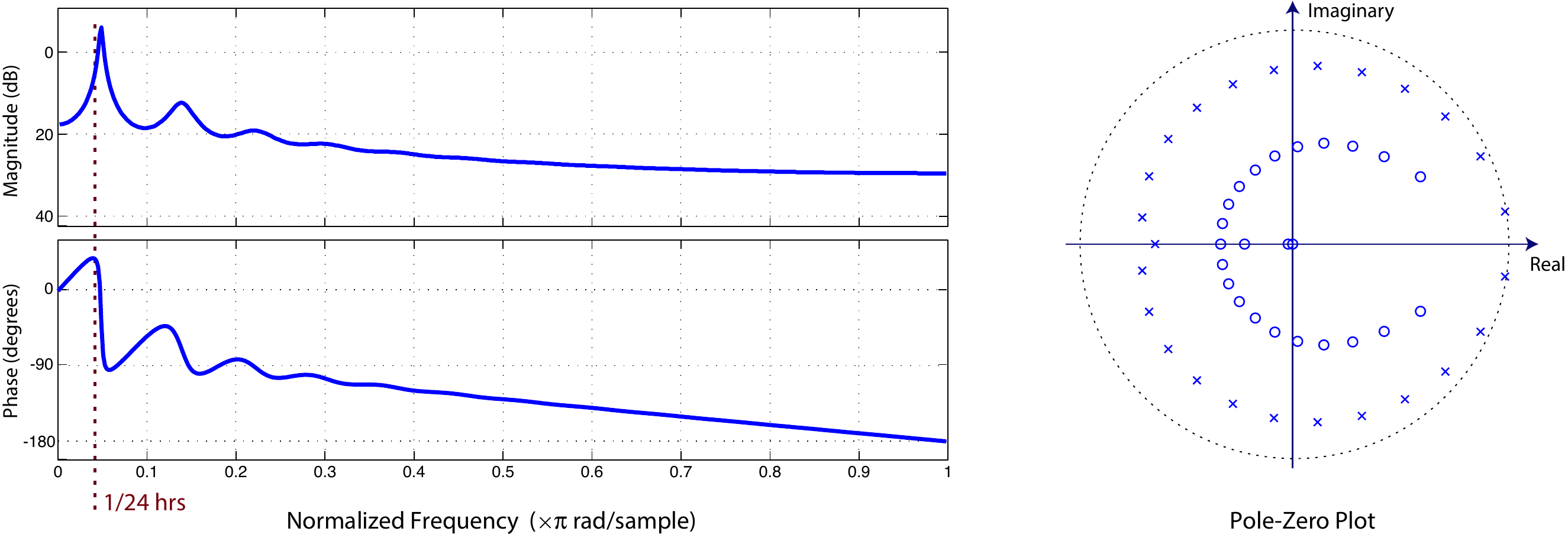}
\vspace{-.05cm} 
\caption{Frequency response and pole-zero plot for the linearized model $C[Iz-A]^{-1}B$}
\vspace{-.25cm}
\label{fig:fr}
\vspace{-08pt}
\end{figure}

Consider now the case in which $\bfzeta$ is a function of time. 
\Fig{fig:fr} shows the   Bode plot and pole-zero plot
for the linear model \eqref{e:LSSmfg}. 
The transfer function from $\bfzeta$ to $\bfgamma$ is BIBO stable and minimum phase.


\notes{I removed the pole-zero cancellations:
The pole-zero cancellations imply that the model cannot be both controllable and observable.
The controllability matrix has rank $23$ and the observability matrix rank $15$.
This is a 50-state model, so the model is neither controllable nor observable.
|}
\notes{See commented text for \it Frequency Sweep (Swept Sine)}

\subsection{Super-sampling}
\label{sec:sim}

Recall the control architecture described at the start of \Section{s:mfg}. 
At any given time, the desired power consumption/curtailment is
determined by the 
BA 
based on its knowledge of dispatachable and uncontrollable generation as well as prediction of load. This is passed through a band-pass filter and scaled appropriately  based on the proportion of ancillary service provided by the pools, and the average power consumption of pool pumps.  The resulting reference signal is denoted $\bfmr$.

We introduce here a refinement of the randomized control scheme to account for delay in the system:  Even if sampling takes place each hour, if a percentage of pools turn off in response to a regulation signal, then the power consumption in the grid will drop nearly instantaneously.  Nevertheless, the control system model will have a one hour delay, which is unacceptable.    

To obtain a more responsive system we employ   ``super-sampling'' at the grid level, which is obtained as follows:
We maintain the assumption  that each pool checks the regulation signal at intervals of length $T$.  
However,  the pools have no common clock.   

It is convenient to model super-sampling via binning of time, so that we retain a discrete time model.
Let  $m>1$ denote a ``super-sampling'' parameter.  At the grid-level the system is in discrete time,  with sampling interval $T/m$.  For example, if $T=30$ minutes, then $m=6$ corresponds to a five minute sampling interval.  A pool is class $i$ if the reference signal is checked at times $nT + (i-1)T/m$,  with $n\ge 0$, $1\le i\le m$.

Letting $y_t^i$ denote the fraction of pools    in the $i$th class that are  operating,  the total that are   operating at time $t$ is  the sum,
\[
y_t =\sum_{i=0}^{m-1} y_t^i
\]
Let $H_0$ denote the discrete time transfer function using $m=1$,
which is simply the transfer function for the linear state space model 
\eqref{e:LSSmfg}.   For general $m$,  the transfer function   from $\zeta$ to $y$ is  
\spm{This needs to be checked - we need a clear analysis in the paper.  }
\begin{equation}
H(z^{-1}) =  z^mH^0(z^{-m}) L(z) 
\label{e:supSampleHvt}
\end{equation}
where $L$ is the low pass filter,
\[
L(z)=  \frac{1}{m} \sum_{i=1}^{m} z^{-i}  =  \frac{1}{m}  z^{-1}  \frac{1-z^{-m}}{1-z^{-1}} 
\]
The term ``$1/m$'' appears because the pools in each bin contribute this fraction of total ancillary service.
In the second representation  there is a pole-zero cancellation at $z=1$.  The filter $L(z)$ has $m-1$ zeros on the unit circle: 
All of the solutions to $z^{m}=1$, except for the solution $z=1$.

Using super-sampling we have achieved our goal of reducing delay:   In   real time, the delay in this model is $T/m$ rather than $T$.

\subsection{Simulation results}

The numerical results described here are based on a stochastic simulation  
of one million pools ($N=10^6$),  using Matlab.   This large number of pools is consistent with Florida or California.  

For the purposes of translation to megawatts,  it is assumed that   each pool in operation consumes $\GNoplus =1$~KW.  Power consumption at time $t$ is assumed to be equal to $N \GNoplus y_t$ (in KW).

\notes{no!! or $10^3 y_t $ (in MW).}

The supersampling approach was used in all of these experiments, with the following values of $T$ and $m$ fixed throughout: Each pool checks the regulation signal every $T=30$ minutes.  The supersampling parameter is  $m=12$, corresponding to $150$~second sampling intervals at the grid level.

The reference signal was chosen to be the BPA regulation signal passed through a low pass filter, shown in \Fig{fig:BPA}.
It was found that one million pools could provide far more regulation than the $\pm$~200~MW required at BPA during this week.  More experiments were conducted in which the signal was scaled to investigate the limits of regulation from a population of one million pools.

We summarize results obtained from
two sets of experiments  conducted in two scenarios.  In the first,  the symmetric model based on the nominal model was used,  with switching probability \eqref{e:pvarrho}.

\begin{figure}[h]
\Ebox{.5}{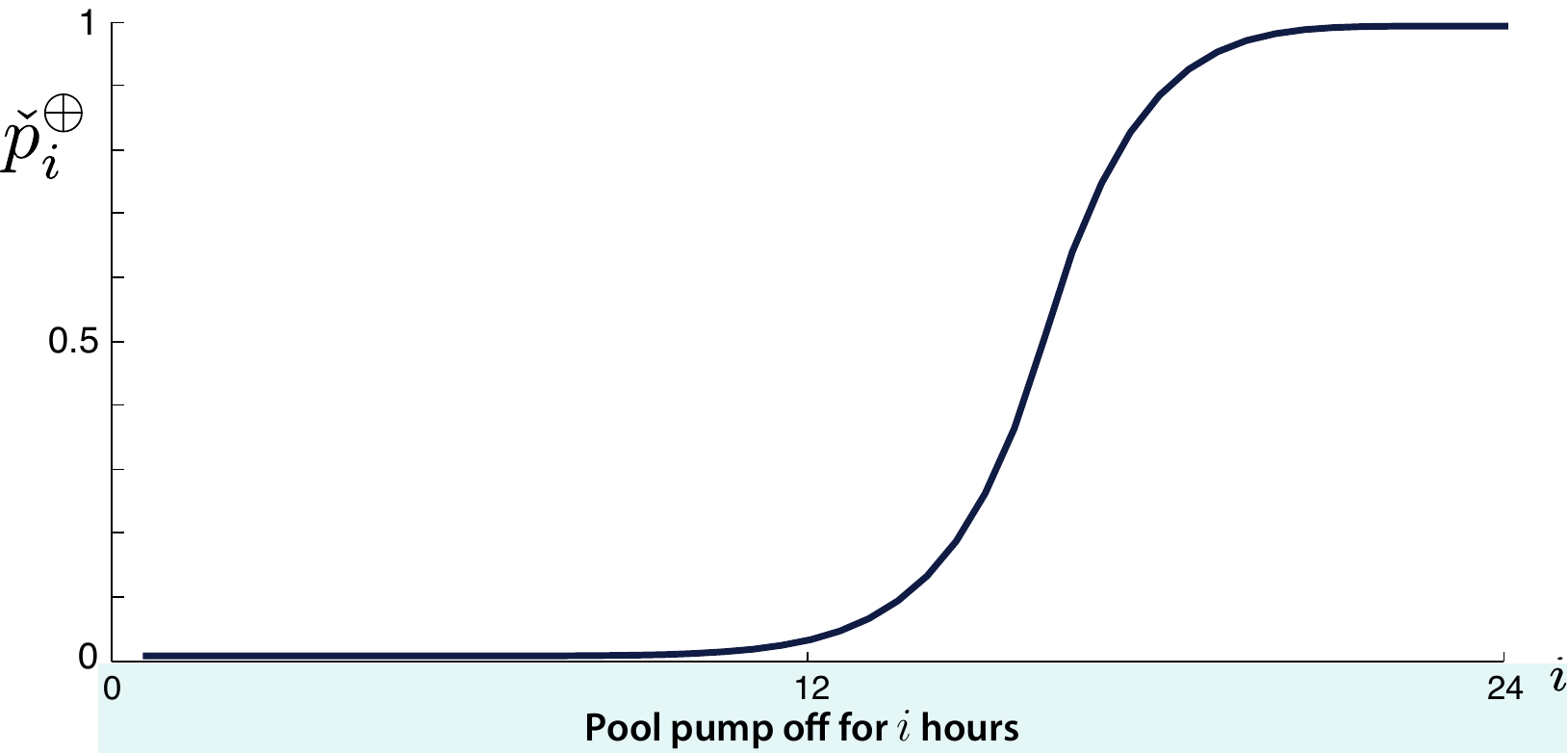}
\vspace{-.05cm}
\caption{Nominal model with 8 hour cleaning schedule.}
\label{fig:plus8}
\vspace{-08pt}
\end{figure}

The second scenario was based on a shorter cleaning schedule of 8 hours per day.   
The switching probabilities 
defined in \eqref{e:pdefn}
were used in which $\alpha=1/3$,  and
$\gamma = 6$ in both scenarios.
The function $ p^\oplus $ using $\alpha=1/3$ is shown in \Fig{fig:plus8}.

In both scenarios, the linearization \eqref{e:LSSmfg} is minimum phase:
All zeros of $H_0(z) = C(Iz-A)^{-1}B$ lie strictly within the unit disk in the complex plane.   With the introduction of super-sampling, the resulting transfer function \eqref{e:supSampleHvt} also has zeros on the unit circle.

In these experiments it was assumed that the BA had perfect measurements of the total power consumption of the population of pools.   PI control was used to obtain the signal $\bfzeta$:  A proportional gain of $20$, and integral gain of $4$ worked well in all cases.   That is, the command $\bfzeta$ was taken to be
\[
\zeta_t = 20 e_t + 4 e^I_t,\qquad \text{with} \ \ e_t = r_t-y_t\ \ \text{\it and} \ \ e^I_t = \sum_{k=0}^t e_k
\]
This is of the form $
\zeta_t = \phi_t(\mu_0,\dots,\mu_t)$,  $ t\ge 0$, that is required in \Theorem{t:MFL}.

 \begin{figure}[h]  
\Ebox{1}{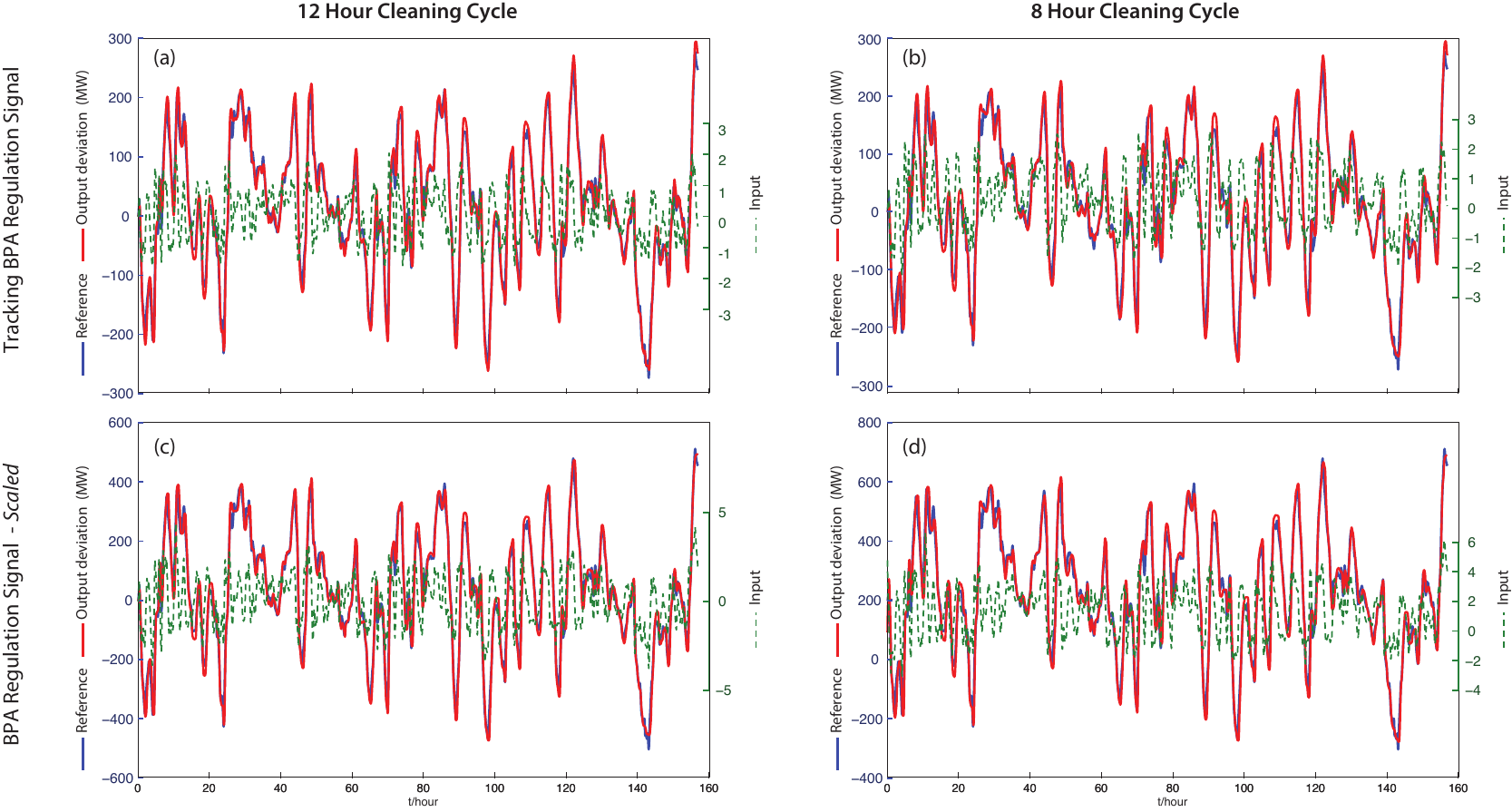}
\vspace{-.25cm} 
\caption{Closed loop simulation in two scenarios, using two different reference signals.  }
\vspace{-.05cm}
\label{f:YueSimScale} 
\end{figure}   

The average proportion of time that a pool is on will be approximately $1/2$ in Scenario 1,  and $1/3$ in Scenario 2.  Consequently, the class of regulation signals that can be tracked is not symmetric  in Scenario 2:  The population of pools has more potential for increasing rather than decreasing power consumption.
 
To attempt to quantify this effect,  define \textit{potential capacity} as the upper and lower limits of power deviation, subject to the constraint that tracking performance does not degrade,
denoted $\{+\text{Demand}, -\text{Supply}\}$. 
Through simulations it was found that 
the potential capacity in Scenario~$1$ is $\{+500MW, -500MW\}$,  and $\{+695MW, -305MW\}$ in  Scenario~$2$.

Results from four experiments are shown in \Fig{f:YueSimScale}.  Subplots (a) and (b) show tracking results using the low-pass filtered signal shown in \Fig{fig:BPA},   and the second row shows tracking performance when the signal magnitude is increased and shifted to match its potential capacity.  The tracking performance is remarkable in all cases.  In particular, it is surprising that a $\pm400$~MW signal can be tracked, given that the average power consumption of the pools is $500$~MW in Scenario~1.

 \begin{figure}[h]  
\Ebox{1}{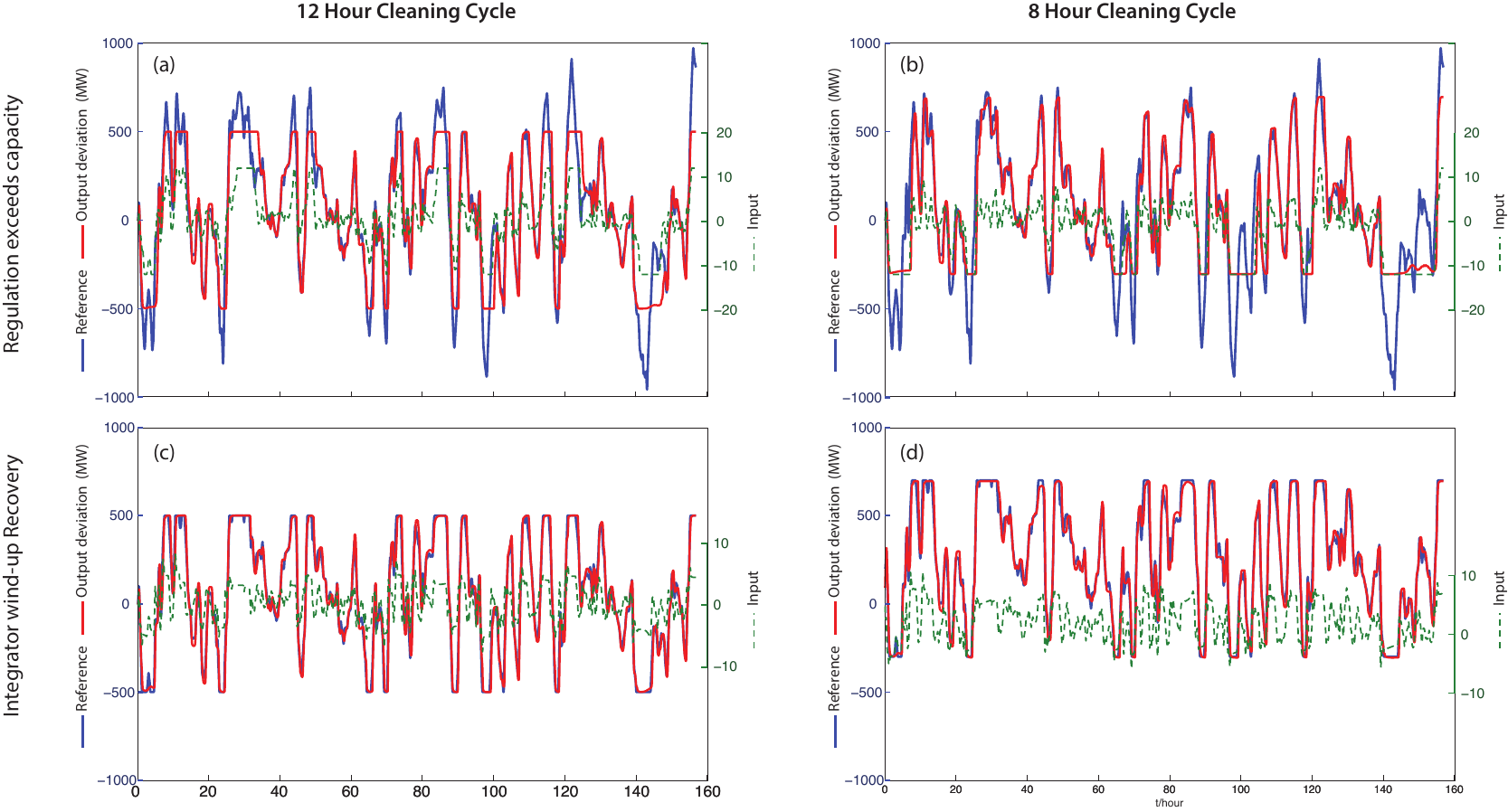}
\vspace{-.25cm} 
\caption{The impact of exceeding capacity}
\vspace{-.05cm}
\label{f:YueSimScaleWindup} 
\end{figure}  


Subplots (a) and (b) in \Fig{f:YueSimScaleWindup} shows what happens when the reference signal exceeds capacity.   Two sources of error are evident in these plots.   First,  the power deviation saturates when all of the $10^6$ pools are turned off, or all are turned on.   Secondly,    large tracking  errors are observed immediately after saturation.  This is a consequence of memory in the PI controller -- what is known as \textit{integrator windup}.   To solve this problem,  the BA should truncate the regulation signal so that it does not exceed the values $\{+\text{Demand}, -\text{Supply}\}$. 
Subplots (c) and (d) in \Fig{f:YueSimScaleWindup} use the same regulation signal used in (a), (b),   but truncated to meet these capacity constraints.  Once again, the tracking is nearly perfect.

\paragraph*{Individual risk}
These simulation experiments have focused on the service to the grid, and the accuracy of the mean-field model approximations.  The fidelity of approximation is remarkable.

The next question to ask is, what happens to an individual pool?  Because of constraints on the regulation signal, it is found in simulations that the average cleaning time for each pool owner is close to the target values (either 12 or 8 hours per day in the two scenarios treated here).   This is to be expected by the Law of Large Numbers.

The Central Limit Theorem can be appealed to if we wish to understand 
the impact of this control architecture  on an individual pool.   In simulations we find that the empirical  distribution of hours cleaned over a four day period appears to be roughly Gaussian, and hence there is a portion of pools that are under-cleaned, and another portion that receive too many hours of cleaning.     

Risk to individual consumers can be reduced or eliminated by using  an additional layer of control at the loads.   If over a period of several days the system detects over or under-cleaning,  then the control system will ignore the signal sent by the BA.   The aggregate impact of this modification represents a small amount of un-modeled dynamics.  In preliminary experiments we have seen virtually no impact on tracking; only a small reduction in capacity.  Analysis of individual risk is a topic of ongoing research.

\section{Conclusions} 
\label{s:conclude}

The simplicity of the MDP solution, and the remarkable accuracy of the LTI approximation for the mean-field model makes this approach appealing for this and many related applications.


There are several issues that have not been addressed here:  
\begin{romannum}
\item 
We do not fully understand the potential cost to consumers in terms of energy, or risk in terms of rare events in which the pool is under- or over-cleaned.  It is likely that hard constraints on performance can be put in place without impacting the analysis.

\item
Does the grid operator need to know the real-time power consumption of the population of pools?  Probably not. The BA is interested in regulating frequency,   and this may be the only measurements needed for harnessing ancillary service from these loads.  The grid frequency passed through a band-pass filter could serve as a surrogate for the measurement $y_t$ assumed in this paper.

It may be valuable to have \textit{two} measurements at each load:  The BA command, and local frequency measurements.

\item
How can we engage consumers?
The formulation of contracts with customers requires a better understanding of the value of ancillary service, as well as consumer preferences.   

\end{romannum}



 


\bibliographystyle{IEEEtran}
\bibliography{strings,markov,q,extras}

\clearpage

 \appendix

\subsection{Controlled transition matrix:  Proof of \Proposition{lem:Pcheck-infinite}.}
 
The dependency of $\eta^*$, $h=\log(v)$, and $\lambda$ on $ \zeta$ will be suppressed to simplify notation in the proof of  \Proposition{lem:Pcheck-infinite} that follows.

The existence of a unique maximal eigenvalue and a positive eigenvector is a consequence of the Perron-Frobenious Theorem \cite{sen81,num84}. Recent results on this and associated multiplicative ergodic theory are contained in \cite{konmey05a}, from which the identity $\Lambda\eqdef\log(\lambda)=\eta^*_\zeta$  can be established (see  \cite[Theorem 1.2]{konmey05a}).    

The bounds in \eqref{e:cpOpt} imply a rate of convergence of the finite horizon cost using $\cP$ to its infinite-horizon limit,  which in particular implies that $\Lambda\eqdef\log(\lambda)=\eta^*_\zeta$.
 To complete the proof of
 \Proposition{lem:Pcheck-infinite} 
 it remains to establish this pair of bounds.
 
We begin with the second bound in \eqref{e:cpOpt}.
Let $\cp^T$ denote the probability on strings on $\state^{T}$ induced by $\cP$, with initial condition $x_0$ given.  This can be expressed,
\[
\cp^T (x_1,\dots, x_{T}) =  \frac{v(x_T)}{v(x_1)}\exp\bigl(\zeta \sum_{t=1}^{T} \util(x_t) - T \Lambda \bigr) p_0 (x_1,\dots, x_{T}) 
\]
Since $\cp^T$ is a probability measure we have,
\[
1 
= 
	\sum_{x_i} \cp^T(x_1,\dots, x_{T})  
=
	\sum_{x_i} \frac{v(x_T)}{v(x_1)}\exp\Bigl(\zeta \sum_{t=1}^{T} \util(x_t) - T \Lambda \Bigr) p_0 (x_1,\dots, x_{T}) 
\]
which gives,
\[
T\Lambda= \log\Bigl\{ \Expect\Bigl[\frac{v(X_T)}{v(X_1)}
\exp\Bigl(\zeta \sum_{t=1}^{T} \util(X_t) \Bigr) \Bigr] \Bigr\}\,.
\]
The bound on $|T\Lambda -  \welf_T(p^{T*})   |$ follows from this identity combined with \Proposition{t:twisted} (which establishes  $\welf_T(p^{T*})=\Lambda_T(\zeta) $).

Substitution of $\cp^T$ into the definition of $\welf_T$ gives,
\[
\welf_T(\cp^T) =  - \Expect_p\Bigl[\log(v(X(T)) -\log(v(X(1)) \Bigr]   + T \Lambda
\]
This combined with the second bound in \eqref{e:cpOpt} gives the desired upper bound on $\welf_T(p^{T*}) -  \welf_T(\cp^{T}) $.
\qed

\subsection{Computation}

The approximations of $\eta^*_\zeta$ and $h^*_\zeta$ are based on Poisson's equation.
In general, we say that \textit{$g$ is the solution to Poisson's equation with forcing function $f$} if the following equation holds:  
\begin{equation}
g- P_0 g  = f
\label{e:fish}
\end{equation}
For this finite state space Markov chain, it follows from invariance of $\pi_0$ that $\pi_0(f) = \sum_x \pi_0(x)f(x)=0$.

We show here that  the second order term $\varH$ defined in \eqref{e:SOvH} can be expressed as a solution to a certain Poisson equation, and explain how to compute solutions.

\subsubsection{Poisson's equation}
\label{s:fish}

This is a finite state Markov chain,  so Perron-Frobenious theory is very simple and attractive  \cite{sen81,num84,konmey05a}.   Let $s\colon\state\to\Re_+$ be a function (not identically zero),
and $\nu$ a probability measure such that the \textit{minorization condition} holds,
\[
   P_0(x,y)\ge s(x)\nu(y),\qquad x,y\in\state
\]
This is written $P\ge s\otimes\nu$.   We will take $s(x)= \ind\{x=\atom\}$ and $\nu(y) = P_0(\atom,y)$,  so that for all $y$,
\[
\barP(x,y)\eqdef
P_0(x,y)-s(x)\nu(y)=\begin{cases}  0 & x=\alpha
\\
P_0(x,y) & x\neq \alpha
\end{cases}
\]

The invariant probability measure $\pi_0$ can be expressed in terms of a matrix inverse.   
We first use the implication,
\[
\pi_0P_0 = \pi_0\Longrightarrow \pi_0(P_0-s\otimes\nu) = \pi_0-\delta \mu
\] 
where $\delta= \sum_x \pi_0(x) s(x) $ is a constant.  We then invert,
\[
\pi_0= \delta  \nu [I-(P_0-s\otimes\nu) ]^{-1} 
\]
We substitute $\barP = P_0 - s\otimes\nu$,  and note that the measure $ \mu = \nu [I-\barP ]^{-1} $ is the unnormalized invariant measure.  Normalization gives the invariant probability measure,
\[
\pi_0(x) = \mu(x)/[\sum_y \mu(y)],\qquad x\in\state.
\]

The inverse  $Z = [I- \barP ]^{-1} $ is called the potential matrix \cite{num84}.

A similar computation is used to solve Poisson's equation.  
Assume that $\nu(\derH) = \sum\nu(x)\derH(x)=0$.  This is without loss of generality.
Under this assumption,
\[
\tilutil +P_0 \derH  = \derH  \Longrightarrow \tilutil +(P_0-s\otimes\nu) \derH  = \derH 
\]
which gives $\derH  = Z \tilutil$.

\subsubsection{Second-order approximation}

Let $h'_\zeta$ and $h''_\zeta$ denote  the first and second derivatives of $h^*_\zeta$ with respect to the argument $\zeta$.  Each of these are real-valued functions on $\state$. 
The Taylor series approximation \eqref{e:happrox}
 in
\Proposition{t:etaSecondOrder} implies that 
\[
h'_\zeta\Big|_{\zeta=0} =\derH \quad \text{\it and}  \quad h''_\zeta\Big|_{\zeta=0} =\varH\,.
\]
Here we obtain representations of these functions that are convenient for computation.  

We show  in 
\Proposition{t:PoissonDerivatives}
that  $\varH$ solves a certain Poisson equation, that can be solved through elementary matrix algebra.
   The  representation for $h''_\zeta$ in 
\Proposition{t:PoissonDerivatives} appears to be new.   

The derivations are most easily obtained using the \textit{nonlinear generator}, defined for any function $g\colon\state\to \Re$ via,
\[
\DV (g) \eqdef   \clP(g) - g  ,  
\]
where $ \clP(g)$  denotes the function,
\[
 \clP(g) \Big|_x \eqdef
\log( P e^g) \Big|_x = \log\Bigl( \sum_y P(x,y) e^{g(y)} \Bigr),\qquad x\in\state
\]
Using this notation, the eigenvector equation can be expressed in a form similar to Poisson's equation,
\begin{equation}
 h^*_\zeta - \clP(h^*_\zeta) = z\util -  \eta^*_\zeta
\label{e:MPE}
\end{equation}
in which $z\util -  \eta^*_\zeta$ plays the role of a forcing function, similar to \eqref{e:fish}.
The similarity between \eqref{e:fish} and \eqref{e:MPE} is why $h^*_\zeta$ is called the solution to the \textit{multiplicative Poisson equation} in \cite{konmey05a}.

We also require a nonlinear operator that defines one-step variance:
For any function $g\colon\state\to \Re$,
\[
\V (g) 
 \Big|_x 
 \eqdef
 \sum_y P_0(x,y) g^2(y) - \Big[  \sum_y P_0(x,y) g(y) \Bigr]^2 
\]
That is, $
\V (g) = P g^2- [Pg]^2
$.

\notes{we should reference
\cite{sch68} -- which considered   a general smooth family of transition matrices. }

\begin{proposition}
\label{t:PoissonDerivatives}
The functions $h'_\zeta$ and $h''_\zeta$ solve the following Poisson equations:
\begin{eqnarray}
\util + \cP_\zeta h'_\zeta 
&=&  h'_\zeta  + \Lambda_\zeta'
\label{e:hPoisson1}
\\
\V(h'_\zeta) +  \cP_\zeta h''_\zeta 
&=&  h''_\zeta  + \Lambda_\zeta''
\label{e:hPoisson2}
\end{eqnarray}
where $\Lambda_\zeta=\log(\lambda_\zeta) = \eta^*_\zeta$.  
The derivatives also satisfy the boundary conditions,
\[
 h'_\zeta (\atom) = 
 h''_\zeta (\atom) =0. 
\]
In particular, with $\zeta=0$ we obtain,
\begin{eqnarray}
\util + P_0 \derH 
&=&  \derH   + \eta_0
\label{e:hPoisson1b}
\\
\V(\derH ) +  P_0 \varH 
&=&  \varH  + \avar
\label{e:hPoisson2b}
\end{eqnarray}
		\qed
\end{proposition}

\IEEEproof 
The boundary conditions hold because $h^*_\zeta(\atom)=\log(v_\zeta(\atom)) = 0$ for all $\zeta$.
This is by construction --- see \eqref{e:PFv} and \eqref{e:PFeta}.

The remaining results are obtained on differentiating each side of  \eqref{e:MPE}. 
 The first derivative of $\clP(h^*_\zeta)$ is obtained from simple calculus:
\begin{equation}
\frac{d}{d\zeta}  \clP(h^*_\zeta)  =
\frac{d}{d\zeta} \log( P e^{h^*_\zeta})  = \frac{1}{P e^{h^*_\zeta}} P \{ e^{h^*_\zeta} h'_\zeta\}  = \cP_\zeta h'_\zeta
\label{e:DVfirstDer}
\end{equation}
On differentiating each side of 
\eqref{e:MPE} we thus obtain
\eqref{e:hPoisson1}.

For the second derivative of $h^*_\zeta$ 
we require the second derivative of $\clP(h^*_\zeta)$.    
A representation will follow from the product rule, using \eqref{e:DVfirstDer}:
\begin{equation}
\frac{d^2}{d\zeta^2}  \clP(h^*_\zeta)   = \cP_\zeta h''_\zeta + \cP_\zeta' h'_\zeta
\label{e:DVsecondDerProd}
\end{equation}
To obtain an expression for the matrix $\cP_\zeta'$, 
observe that for any $x,y\in\state$ for which $P_0(x,y)> 0$,
\[
\begin{aligned}
\frac{1}{ \cP_\zeta(x,y) }
\frac{d}{d\zeta}  \cP_\zeta(x,y) 
&= 
\frac{d}{d\zeta} \log\bigl( \cP_\zeta(x,y) \bigr) 
\\
&= \frac{d}{d\zeta} \log\bigl(  P(x,y) e^{h^*_\zeta(y) } / Pe^{h^*_\zeta }\, (x) \bigr)
\\
  &=  h'_\zeta(y) -  \frac{d}{d\zeta} \clP(h^*_\zeta) \Big|_x  =  h'_\zeta(y) - \cP_\zeta h'_\zeta\,(x)
\end{aligned}
\]
where the final equality follows from \eqref{e:DVfirstDer}.
This gives,
\[
\frac{d}{d\zeta}  \cP_\zeta(x,y) 
=  \cP_\zeta(x,y) [h'_\zeta(y) - \cP_\zeta h'_\zeta\,(x)],\qquad x,y\in\state.
\] 
Substituting this into \eqref{e:DVsecondDerProd}
gives,
\[
\frac{d^2}{d\zeta^2}  \clP(h^*_\zeta) \Big|_x
=
 \cP_\zeta h''_\zeta \,(x) + \cP_\zeta' h'_\zeta\,(x)
 = 
 \cP_\zeta h''_\zeta \,(x) + \sum_y  \cP_\zeta(x,y) [h'_\zeta(y) - \cP_\zeta h'_\zeta\,(x)]h'_\zeta\,(y)
\]
which simplifies to,
\begin{equation}
\frac{d^2}{d\zeta^2}  \clP(h^*_\zeta)   = \cP_\zeta h''_\zeta + \V (h'_\zeta) 
\label{e:DVsecondDer}
\end{equation}
Differentiating   each side of  \eqref{e:MPE} twice thus gives
\eqref{e:hPoisson2}.
\qed

\end{document}